\pdfoutput=1

\documentclass[letterpaper, 11pt]{article}
\usepackage{amssymb,amsmath}
\usepackage{epsf}
\usepackage{times,bm,mathrsfs}  
\usepackage{amsmath}
\usepackage{amssymb}
\usepackage{amsfonts}
\usepackage{mathrsfs}
\usepackage{bbm}
\usepackage{color}
\usepackage{graphicx}
\usepackage{amscd}
\usepackage{epsfig}

\newcommand{\itemname}[1]{$\mathrm{[#1]}$}
\newcommand{\onetwo}{\iota} 

\newcommand{\quartet}{\qcal}

\newcommand{\forestno}{\mathcal{F}}
\newcommand{\metricno}{\mathcal{D}}

\newcommand{\distmet}{\textsc{DistortedMetric}}

\definecolor{Red}{rgb}{1,0,0}
\definecolor{Blue}{rgb}{0,0,1}
\definecolor{Olive}{rgb}{0.41,0.55,0.13}
\definecolor{Green}{rgb}{0,1,0}
\definecolor{MGreen}{rgb}{0,0.8,0}
\definecolor{DGreen}{rgb}{0,0.55,0}
\definecolor{Yellow}{rgb}{1,1,0}
\definecolor{Cyan}{rgb}{0,1,1}
\definecolor{Magenta}{rgb}{1,0,1}
\definecolor{Orange}{rgb}{1,.5,0}
\definecolor{Violet}{rgb}{.5,0,.5}
\definecolor{Purple}{rgb}{.75,0,.25}
\definecolor{Brown}{rgb}{.75,.5,.25}
\definecolor{Grey}{rgb}{.5,.5,.5}
\definecolor{Black}{rgb}{0,0,0}

\newcommand{\acal}{\mathcal{A}}

\newcommand{\dcal}{\mathcal{D}}
\newcommand{\ecal}{\mathcal{E}}
\newcommand{\fcal}{\mathcal{F}}

\newcommand{\jcal}{\mathcal{J}}
\newcommand{\lcal}{\mathcal{L}}

\newcommand{\qcal}{\mathcal{Q}}
\newcommand{\rcal}{\mathcal{R}}
\newcommand{\scal}{\mathcal{S}}
\newcommand{\tcal}{\mathcal{T}}

\newcommand{\vcal}{\mathcal{V}}

\newcommand{\zcal}{\mathcal{Z}}

\newcommand{\real}{\mathbb{R}}

\newcommand{\eps}{\varepsilon}

\newcommand{\ind}{\mathbbm{1}}

\newcommand{\bdm}{\begin{displaymath}}
\newcommand{\edm}{\end{displaymath}}
\newcommand{\bea}{\begin{eqnarray*}}
\newcommand{\eea}{\end{eqnarray*}}
\newcommand{\bean}{\begin{eqnarray}}
\newcommand{\eean}{\end{eqnarray}}

\newcommand{\prob}{\mathbb{P}}
\newcommand{\expec}{\mathbb{E}}
\newcommand{\var}{\mathrm{Var}}

\newcommand{\bfmu}{\boldsymbol{\mu}}

\newcommand{\poly}{\mathrm{poly}}

\newtheorem{theorem}{Theorem}
\newtheorem{proposition}{Proposition}

\newtheorem{definition}{Definition}
\newtheorem{example}{Example}
\newtheorem{lemma}{Lemma}

\newcommand{\ignore}[1]{}

\renewcommand{\poly}{{\mbox{{\rm poly}}}}

 \def\eps{\varepsilon}

\def\R{\hbox{I\kern-.2em\hbox{R}}}

\def\|{\, | \, }
\def\v0{{\bf 0}}

\def\0{\hat{0}}
\def\1{\hat{1}}

\def\phi{\varphi}

\def\be{\begin{equation}}
\def\ee{\end{equation}}

\def\part{{\cal P}}

\def\R{\mathcal{R}}

\def\eps{\varepsilon}

\newcommand{\weight}{\tau}
\newcommand{\eweight}{\hat\tau}
\newcommand{\path}{\mathrm{Path}}
\newcommand{\dist}{\tau}

\newcommand{\phy}{\tcal}
\newcommand{\rt}{\rho}
\newcommand{\hmgt}[1]{T^{(#1)}}
\newcommand{\hmgv}[1]{V^{(#1)}}
\newcommand{\hmgl}[1]{L^{(#1)}}
\newcommand{\hmgll}[2]{L^{(#1)}_{#2}}
\newcommand{\hmge}[1]{E^{(#1)}}

\newcommand{\hmgphy}[1]{\phy^{(#1)}}
\newcommand{\hmgrt}[1]{\rt^{(#1)}}
\newcommand{\states}{\Phi}
\newcommand{\nstates}{\phi}
\newcommand{\trees}{\mathbb{T}}

\newcommand{\s}{\sigma}

\newcommand{\evr}{\nu}
\newcommand{\rates}{\mathbb{Q}}
\newcommand{\gtr}{\mathbb{G}}
\newcommand{\law}{\lcal}
\newcommand{\sphy}{\mathbb{Y}}
\newcommand{\shmgphy}{\mathbb{HY}}

\newcommand{\corr}{\widehat F}
\newcommand{\estim}{\dcal}
\newcommand{\deep}{\overline{\mathbb{FP}}}
\newcommand{\estdist}{\overline{\mathrm{d}}}
\newcommand{\longtest}{\overline{\mathbb{SD}}}
\newcommand{\triplet}{\mathbb{O}}
\newcommand{\ball}{\mathcal{B}}

\newcommand{\unit}{\bar{e}}

\newcommand{\loglog}[1]{\lfloor #1 \log_2\log_2 n \rfloor}
\newcommand{\LogLog}[1]{\log_2^{#1} n}

\newcommand{\quantum}{\Delta}

\newenvironment{app-proof}[1]{\noindent{\textbf{Proof of #1:}}}{$\blacksquare$\vskip\belowdisplayskip}

\author{
S\'ebastien Roch\footnote{
Department of Mathematics, UCLA.
}
}

\title{\vspace{0cm}
Sequence-Length Requirement of Distance-Based Phylogeny Reconstruction:
Breaking the Polynomial Barrier\footnote{The current manuscript is the full 
version with proofs of~\cite{Roch:08}.
In subsequent work~\cite{Roch:09} the results stated here 
were improved to logarithmic sequence length, thereby matching the best results
for general methods.}
}

\begin{document}

\maketitle

\thispagestyle{empty}

\begin{abstract}
We introduce a new distance-based phylogeny reconstruction technique 
which provably achieves, at sufficiently short branch lengths,
a polylogarithmic sequence-length 
requirement---improving significantly over previous polynomial 
bounds for distance-based methods. 
The technique is based on an averaging procedure that implicitly 
reconstructs ancestral sequences.

In the same token, we extend previous results on phase transitions in phylogeny reconstruction to general
time-reversible models. More precisely, we show that in the so-called Kesten-Stigum zone (roughly,
a region of the parameter space where ancestral sequences are well approximated by ``linear combinations''
of the observed sequences) sequences of length $\poly(\log n)$ suffice for reconstruction
when branch lengths are discretized.
Here $n$ is the number of extant species. 

Our results challenge, to some extent, the conventional wisdom that estimates of evolutionary distances alone 
carry significantly
less information about phylogenies than full sequence datasets.
\end{abstract}


\clearpage

\section{Introduction}

The evolutionary history of a group of organisms is generally represented
by a {\em phylogenetic tree} or {\em phylogeny}~\cite{Felsenstein:04, SempleSteel:03}.
The leaves of the tree represent the current species. 
Each branching indicates a speciation event.
Many of the most popular techniques for reconstructing phylogenies from
molecular data, e.g.~UPGMA, Neighbor-Joining, and BIO-NJ~\cite{SokalSneath:63,
SaitouNei:87, Gascuel:97}, 
are examples of what are known as {\em distance-matrix methods}.
The main advantage of these methods is their speed, 
which stems from a straightforward approach: 
1) the estimation of a {\em distance matrix} from observed molecular sequences;
and 2) the repeated agglomeration of the closest clusters of species.
Each entry of the distance matrix 
is an estimate of the evolutionary distance between the corresponding
pair of species, that is, roughly the time elapsed since their most recent common ancestor.
This estimate is typically obtained by comparing aligned homologous DNA sequences extracted
from the extant species---the basic insight being, the closer the species, the more
similar their sequences.
Most distance methods run in time polynomial in $n$, the number
of leaves, and in $k$, the sequence length. This performance compares very favorably to 
that of the other two main classes of reconstruction methods, likelihood and parsimony methods, 
which are known to be computationally 
intractable~\cite{GrahamFoulds:82, DaySankoff:86, Day:87, MosselVigoda:05,
ChorTuller:06, Roch:06}. 

The question we address in this paper is the following: Is there a price to pay
for this speed and simplicity? There are strong combinatorial~\cite{StHePe:88} 
and statistical~\cite{Felsenstein:04} reasons
to believe that distance methods are not as accurate as more elaborate reconstruction techniques,
notably maximum likelihood estimation (MLE). 
Indeed, 
in a typical instance of the phylogenetic reconstruction problem, we are given
{\em aligned DNA sequences} $\{(\xi^i_l)_{i=1}^k\}_{l\in L}$, one sequence for each leaf $l \in L$,
from which we seek to infer the phylogeny on $L$.
Generally, all {\em sites} $(\xi^i_l)_{l\in L}$, for $i=1,\ldots,k$, are assumed to be
independent and identically distributed according to a Markov model on a tree 
(see Section~\ref{section:definitions}). 
For a subset $W \subseteq L$, we denote by
$\mu_W$ the distribution of $(\xi^i_l)_{l\in W}$ under this model.
Through their use of the distance matrix, distance methods
reduce the data to {\em pairwise 
sequence correlations}, that is, they only use estimates
of $\bfmu_2 = \{\mu_W\ :\ W\subseteq L,\ |W| = 2\}$.
In doing so, they  seemingly fail to take into account more subtle patterns
in the data involving three or more species at a time.
In contrast, MLE for example outputs a model 
that maximizes the {\em joint probability of all observed sequences}. 
We call methods that explicitly use the full dataset, such as MLE, {\em holistic methods}.

It is important to note that the issue is not one of {\em consistency}:
when the sequence length tends to infinity,
the estimate provided by distance methods---just like MLE---typically converges to the correct phylogeny.
In particular, under mild assumptions, it suffices to know the pairwise site distributions $\bfmu_2$
to recover the topology of the phylogeny~\cite{ChangHartigan:91,Chang:96}.
Rather the question is: how fast is this convergence?
Or more precisely, how should $k$ scale as a function of $n$ to guarantee a correct
reconstruction with high probability?
And are distance methods significantly slower to converge than holistic methods?
Although we do not give a complete answer to these questions of practical interest here,
we do provide strong evidence that some of the suspicions against distance methods 
are based on a simplistic view of the distance matrix. In particular,
we open up the surprising possibility that distance methods actually exhibit
optimal convergence rates.

\paragraph{Context.}
It is well-known that some of the most popular distance-matrix methods actually suffer from a
prohibitive sequence-length requirement~\cite{Atteson:99, LaceyChang:06}.
Nevertheless, over the past decade, 
much progress has been made in the design of fast-converging 
distance-matrix techniques, starting with the seminal
work of Erd\"os et al.~\cite{ErStSzWa:99a}. 
The key insight behind 
the algorithm in~\cite{ErStSzWa:99a}, often dubbed 
the Short Quartet Method (SQM), 
is that it discards long evolutionary distances, which are known to be statistically unreliable. 
The algorithm works by first building subtrees 
of small diameter and, in a second stage, putting the pieces back together. 
The SQM algorithm runs in polynomial time and guarantees the correct reconstruction
with high probability
of any phylogeny (modulo reasonable assumptions) when $k = \poly(n)$.
This is currently the best known convergence rate for distance methods.
(See also~\cite{DaMoRo:06,DHJMMR:06,Mossel:07,GrMoSn:08,DaMoRo:08a} 
for faster-converging algorithms involving {\em partial} reconstruction
of the phylogeny.)

Although little is known about the sequence-length requirement of 
MLE~\cite{SteelSzekely:99, SteelSzekely:02},
recent results of Mossel~\cite{Mossel:04a},
Daskalakis et al.~\cite{DaMoRo:06,DaMoRo:08b}, 
and Mihaescu et al.~\cite{MiHiRa:09} on a conjecture of Steel~\cite{Steel:01}
indicate that convergence rates as low as $k=O(\log n)$ can be achieved
when the branch lengths are sufficiently short and discretized, 
using insights from statistical physics.
We briefly describe these results. 

As mentioned above, the classical model of DNA sequence evolution is a Markov model on a tree
that is closely related to stochastic models used to study 
particle systems~\cite{Liggett:85,Georgii:88}. This type of model undergoes
a phase transition that has been extensively studied in 
probability theory and statistical physics: 
at short branch lengths (in the binary symmetric case, up to
15\% divergence {\em per edge}), in what is called the {\em reconstruction
phase}, good estimates of the ancestral sequences can be obtained from the
observed sequences; on the other hand, outside the reconstruction phase,
very little information about ancestral states diffuses to the leaves.
See e.g.~\cite{EvKePeSc:00} and references therein.
The new algorithms in~\cite{Mossel:04a,DaMoRo:06,DaMoRo:08b,MiHiRa:09}
exploit this phenomenon by alternately 1) reconstructing a few levels of the tree
using distance-matrix techniques and 2) estimating distances between {\em internal} nodes
by reconstructing ancestral sequences at the newly uncovered nodes. 
The overall algorithm is \emph{not} distance-based, however, as the ancestral
sequence reconstruction is performed using a complex function
of the observed sequences named {\em recursive majority}.
The rate $k = O(\log n)$ achieved by these algorithms is known
to be necessary in general. Moreover, the slower rate $k=\poly(n)$ is in fact 
necessary for all methods---distance-based or holistic---outside 
the reconstruction phase~\cite{Mossel:03}.
In particular, note that distance methods are in some sense ``optimal'' {\em outside} the
reconstruction phase by the results of~\cite{ErStSzWa:99a}.

\paragraph{Beyond the oracle view of the distance matrix.}
It is an outstanding open problem to determine whether distance methods
can achieve $k = O(\log n)$ in the reconstruction phase\footnote{
Mike Steel offers a 100\$ reward for the solution of this problem. 
}.
From previous work on fast-converging distance methods, it is tempting
to conjecture that $k = \poly(n)$ is the best one can hope for.
Indeed, all previous algorithms use the following
``oracle view'' of the distance matrix, as formalized by King et al.~\cite{KiZhZh:03}
and Mossel~\cite{Mossel:07}.
As mentioned above, the reliability of distance estimates depends on the true
evolutionary distances. From standard concentration inequalities,
it follows that
if leaves $a$ and $b$ are at distance $\weight(a,b)$, then the usual distance
estimate $\eweight(a,b)$ (see Section~\ref{section:definitions}) 
satisfies:
\begin{equation}\label{eq:oracle}
\text{if}\ \weight(a,b) < D + \eps\ \text{or}\ \eweight(a,b) < D + \eps
\ \text{then}\ |\weight(a,b) - \eweight(a,b)| < \eps,
\end{equation} 
for $\eps, D$ such that
$k \propto (1 - e^{-\eps})^{-2} e^{2D}$.
Fix $\eps > 0$ small and $k \ll \poly(n)$. 
Let $T$ be a complete binary tree with $\log_2 n$ levels. Imagine that the distance matrix
is given by the following oracle: on input a pair of leaves $(a,b)$ the oracle returns an estimate 
$\eweight(a,b)$ which satisfies (\ref{eq:oracle}). Now, notice that
for any tree $T'$ which is identical to $T$ on the first $\log_2 n/2$ levels
above the leaves, the oracle is allowed to return the same distance estimate as for $T$.
That is, we cannot distinguish $T$ and $T'$ in this model unless $k  = \poly(n)$.
(This argument can be made more formal along the lines of~\cite{KiZhZh:03}.)

What the oracle model ignores is that, under the assumption that the sequences are generated
by a Markov model of evolution, the distance estimates $(\eweight(a,b))_{a,b\in [n]}$
are in fact {\em correlated random variables}. More concretely, for leaves
$a$, $b$, $c$, $d$, note that the joint distribution of $(\eweight(a,b), \eweight(c,d))$
depends in a nontrivial way on the joint site distribution 
$\mu_W$
at $W = \{a,b,c,d\}$.
In other words,
even though the distance matrix is---seemingly---only an estimate of the pairwise
correlations $\bfmu_2$, it actually contains {\em some} information about all joint distributions.   
Note however that it is not immediately clear how to exploit this extra information or even how useful it could be.

As it turns out, the correlation structure of the distance matrix 
is in fact \emph{very informative} at short branch lengths.
More precisely, we introduce in this paper a new distance-based method 
with a convergence rate
of $k = \poly(\log n)$ 
in the reconstruction phase (to be more accurate,
in the so-called Kesten-Stigum phase; see below)---improving significantly over
previous $\poly(n)$ results. 
Note that the oracle model allows only the reconstruction
of a $o(1)$ fraction of the levels in that case. 
Our new algorithm involves a distance averaging
procedure that implicitly reconstructs ancestral sequences, thereby taking advantage
of the phase transition discussed above. 
We also obtain the first results on
Steel's conjecture beyond the simple symmetric models
studied by Daskalakis et al.~\cite{DaMoRo:06,DaMoRo:08b,MiHiRa:09} 
(the so-called CFN and Jukes-Cantor models). 
In the next subsections, we introduce general definitions and state our results
more formally. We also give an overview of the proof.

\paragraph{Further related work.}
For further related work on efficient phylogenetic tree reconstruction, 
see~\cite{ErStSzWa:99b, HuNeWa:99, CsurosKao:01, Csuros:02}.

\subsection{Definitions}\label{section:definitions}

\paragraph{Phylogenies.} We define phylogenies and evolutionary distances more formally.
\begin{definition}[Phylogeny]
A {\em phylogeny} is a rooted, edge-weighted, leaf-labeled tree 
$\phy = (V,E,[n],\rt;\weight)$
where:
$V$ is the set of vertices;
$E$ is the set of edges;
$L = [n] = \{0,\ldots,n-1\}$ is the set of leaves;
$\rt$ is the root;
$\weight : E \to (0,+\infty)$ is a positive edge weight function.
We further assume that all internal nodes in $\phy$ have degree $3$
except for the root $\rt$ which has degree $2$. 
We let $\sphy_n$ be the set of all such phylogenies on $n$ leaves
and we denote $\sphy = \{\sphy_n\}_{n\geq 1}$.
\end{definition}
\begin{definition}[Tree Metric]
For two leaves $a,b \in [n]$, we denote by $\path(a,b)$
the set of edges on the unique path between $a$ and $b$.
A {\em tree metric} on a set $[n]$ is a positive function
$d:[n]\times[n] \to (0,+\infty)$ such that there exists 
a tree $T = (V,E)$ with leaf set $[n]$ and an edge weight 
function $w:E \to (0,+\infty)$ satisfying the following: for all leaves
$a,b \in [n]$
\begin{equation*}
d(a,b) = \sum_{e\in \path(a,b)} w_e.
\end{equation*}
For convenience, we denote by $\left(\dist(a,b)\right)_{a,b\in [n]}$
the tree metric corresponding to the phylogeny $\phy = (V,E,[n],\rt;\weight)$.
We extend $\weight(u,v)$ to all vertices $u,v \in V$ in the obvious way.
\end{definition}
\begin{example}[Homogeneous Tree]\label{ex:homo}
For an integer $h \geq 0$, we denote by 
$\hmgphy{h} = (\hmgv{h}, \hmge{h}, \hmgl{h}, \hmgrt{h}; \weight)$
a rooted phylogeny where $\hmgt{h}$ is the $h$-level
complete binary tree with arbitrary edge weight function
$\weight$ and $\hmgl{h} = [2^h]$. 
For $0\leq h'\leq h$, we let $\hmgll{h}{h'}$ be the
vertices on level $h - h'$ (from the root). In particular,
$\hmgll{h}{0} = \hmgl{h}$ and $\hmgll{h}{h} = \{\hmgrt{h}\}$.
We let $\shmgphy = \{\shmgphy_n\}_{n\geq 1}$ be the set of all phylogenies
with homogeneous underlying trees.
\end{example}

\paragraph{Model of molecular sequence evolution.}
Phylogenies are reconstructed from molecular sequences
extracted from the observed species. The standard model
of evolution for such sequences is a Markov model
on a tree (MMT). 
\begin{definition}[Markov Model on a Tree]
Let $\states$ be a finite set of character states with $\nstates = |\states|$. 
Typically $\states = \{+1,-1\}$ or $\states = \{\mathrm{A}, \mathrm{G},
\mathrm{C}, \mathrm{T}\}$. Let $n \geq 1$ and let $T = (V,E,[n],\rt)$
be a rooted tree with leaves labeled in $[n]$. 
For each edge $e \in E$, we are given a $\nstates\times\nstates$
stochastic matrix 
$M^e = (M^e_{ij})_{i,j \in \states}$,
with fixed stationary distribution 
$\pi = (\pi_i)_{i\in \states}$.
An MMT $(\{M^e\}_{e\in E}, T)$ 
associates a state $\s_v$ in $\states$ to each
vertex $v$ in $V$ as follows: 
pick a state for the root $\rt$ according to $\pi$;
moving away from the root, choose a state for each 
vertex $v$ independently according to the distribution 
$(M^e_{\s_u, j})_{j\in\states}$,
with $e = (u,v)$ where $u$ is the parent of $v$. 
\end{definition}
The most common MMT used in phylogenetics is the
so-called general time-reversible (GTR) model. 
\begin{definition}[GTR Model]\label{def:gtr}
Let $\states$ be a set of character states with $\nstates = |\states|$
and $\pi$ be a distribution on $\states$ satisfying $\pi_i > 0$ for all
$i\in\states$.  
For $n \geq 1$, let $\phy = (V,E,[n],\rt;\weight)$ be a phylogeny.
Let $Q$ be a $\nstates\times\nstates$ rate matrix, that is,
$Q_{ij} > 0$ for all $i\neq j$ and
$$\sum_{j\in \states} Q_{ij} = 0,$$
for all $i \in \states$.
Assume $Q$ is reversible with respect to $\pi$, that is,
$$\pi_i Q_{ij} = \pi_j Q_{ji},$$
for all $i,j \in \states$.
The GTR model on $\phy$ with rate matrix $Q$ is an MMT 
on $T = (V,E,[n], \rt)$ with transition matrices
$M^e = e^{\weight_e Q}$,
for all $e\in E$.
By the reversibility assumption, $Q$ has $\nstates$ 
real eigenvalues
$$0 = \Lambda_1 > \Lambda_2 \geq \cdots \geq \Lambda_{\nstates}.$$
We normalize $Q$ by fixing $\Lambda_2 = -1$.
We denote by $\rates_\nstates$ the set of all such rate matrices.
We let $\gtr_{n,\nstates} = \sphy_n \otimes \rates_\nstates$ 
be the set of all $\nstates$-state GTR models on 
$n$ leaves. We denote $\gtr_\nstates = 
\left\{\gtr_{n,\nstates}\right\}_{n \geq 1}$.
We denote by $\xi_W$ the vector of states on the vertices
$W\subseteq V$. In particular, $\xi_{[n]}$ are the states
at the leaves. 
We denote by $\law_{\phy,Q}$ the distribution of $\xi_{[n]}$.
\end{definition}
GTR models include as special cases many popular models such as the CFN model.
\begin{example}[CFN Model]\label{ex:cfn}
The {\em CFN model} is the GTR model with $\nstates = 2$,
$\pi = (1/2, 1/2)$,
and
\begin{equation*}
Q 
=
Q^{\mathrm{CFN}}
\equiv
\left(
\begin{array}{cc}
-1/2 & 1/2\\
1/2 & -1/2
\end{array}
\right).
\end{equation*}
\end{example}
\begin{example}[Binary Asymmetric Channel]
More generally, letting $\states = \{+,-\}$ and
$\pi = (\pi_{+}, \pi_{-})$,
with $\pi_{+},\pi_{-} > 0$, we can take
\begin{equation*}
Q 
=
\left(
\begin{array}{cc}
-\pi_{-} & \pi_{-}\\
\pi_{+} & -\pi_{+}
\end{array}
\right).
\end{equation*}
\end{example}

\paragraph{Phylogenetic reconstruction.}
A standard assumption in molecular evolution
is that each site in a sequence (DNA, protein, etc.)
evolves {\em independently} according to a
Markov model on a tree, such as the GTR model above.
Because of the reversibility assumption, the root
of the phylogeny cannot be identified and we reconstruct phylogenies
up to their root.
\begin{definition}[Phylogenetic Reconstruction Problem]
Let $\widetilde\sphy = \{\widetilde\sphy_n\}_{n\geq 1}$ 
be a subset of phylogenies and 
$\widetilde\rates_\nstates$ be a subset of rate matrices on $\nstates$ states. 
Let $\phy = (V,E,[n],\rt;\weight) \in \widetilde\sphy$. 
If $T = (V,E,[n],\rt)$ is the rooted tree underlying $\phy$,
we denote by $T_{-}[\phy]$ the tree $T$ where the root
is removed: that is, we replace the two edges adjacent to the root
by a single edge. We denote by $\trees_n$ the set of all
leaf-labeled trees on $n$ leaves with internal degrees $3$
and we let $\trees = \{\trees_n\}_{n\geq 1}$. 
A {\em phylogenetic
reconstruction algorithm} is a collection of maps 
$\acal = \{\acal_{n,k}\}_{n,k \geq 1}$ 
from sequences $(\xi^i_{[n]})_{i=1}^k \in (\states^{[n]})^k$ 
to leaf-labeled trees $T \in \trees_n$.
We only consider algorithms $\acal$ computable in time polynomial
in $n$ and $k$.
Let $k(n)$ be an increasing function of $n$. We say that $\acal$
solves the {\em phylogenetic reconstruction problem} 
on $\widetilde\sphy \otimes \widetilde\rates_\nstates$ with sequence length $k = k(n)$
if for all
$\delta > 0$, there is $n_0 \geq 1$ such that for all $n \geq n_0$,
$\phy \in \widetilde\sphy_n$, $Q \in \widetilde\rates_\nstates$,  
\begin{equation*}
\prob\left[\acal_{n,k(n)}\left((\xi^i_{[n]})_{i=1}^{k(n)}\right) = 
T_-[\phy]\right] \geq 1 - \delta,
\end{equation*}
where $(\xi^i_{[n]})_{i=1}^{k(n)}$ are i.i.d.~samples from $\law_{\phy,Q}$. 
\end{definition} 
An important result of this kind was given by Erdos et al.~\cite{ErStSzWa:99a}.
\begin{theorem}[Polynomial Reconstruction~\cite{ErStSzWa:99a}]\label{thm:essw}
Let $0 < f \leq g < +\infty$ and denote by $\sphy^{f,g}$ the set of 
all phylogenies $\phy = (V,E,[n],\rt;\weight)$ satisfying
$f \leq \weight_e \leq g,\ \forall e\in E$.
Then, for all $\nstates \geq 2$ and all $0 < f \leq g < +\infty$,
the phylogenetic reconstruction problem on $\sphy^{f,g}\otimes\rates_\nstates$
can be solved with $k = \poly(n)$.
\end{theorem}
This result was recently improved by 
Daskalakis et al.~\cite{DaMoRo:06,DaMoRo:08b} (see also~\cite{MiHiRa:09})
in the so-called Kesten-Stigum reconstruction phase, that is, when
$g < \ln\sqrt{2}$.
\begin{definition}[$\quantum$-Branch Model]
Let $0 < \quantum \leq f \leq g < +\infty$ and denote by $\sphy^{f,g}_\quantum$ the set of 
all phylogenies $\phy = (V,E,[n],\rt;\weight)$ satisfying
$f \leq \weight_e \leq g$ where $\weight_e$ is an integer multiple of 
$\quantum$, for all $e\in E$.
For $\nstates \geq 2$ and $Q\in \rates_{\nstates}$,
we call $\sphy^{f,g}_\quantum \otimes \{Q\}$ 
the $\quantum$-Branch Model ($\quantum$-BM).
\end{definition}
\begin{theorem}[Logarithmic Reconstruction~\cite{DaMoRo:06,DaMoRo:08b}. See also~\cite{MiHiRa:09}.]\label{thm:opt}
Let $g^* = \ln \sqrt{2}$.
Then, for all $0 < \quantum \leq f \leq g < g^*$,
the phylogenetic reconstruction problem on $\sphy^{f,g}_\quantum \otimes \{Q^{\mathrm{CFN}}\}$
can be solved with $k = O(\log n)$\footnote{
The correct statement of this result appears in~\cite{DaMoRo:08b}.
Because of different conventions, our edge weights are scaled by a factor of
$2$ compared to those in~\cite{DaMoRo:08b}. The dependence of $k$ in $\quantum$
is $\quantum^{-2}$.
}.
\end{theorem}

\paragraph{Distance methods.} The proof of Theorem~\ref{thm:essw} uses
{\em distance methods}, which we now define formally. 
\begin{definition}[Correlation Matrix]\label{def:distest}
Let $\states$ be a finite set with $\nstates \geq 2$.
Let 
$(\xi_a^i)_{i=1}^k, (\xi_b^i)_{i=1}^k \in \states^k$
be the sequences at $a, b \in [n]$.
For $\upsilon_1, \upsilon_2 \in \states$, we define the 
{\em correlation matrix} between $a$ and $b$ by
\begin{equation*}
\corr^{ab}_{\upsilon_1 \upsilon_2} 
= \frac{1}{k} \sum_{i=1}^k \ind\{\xi_a^i = \upsilon_1, \xi_b^i = \upsilon_2\},
\end{equation*}
and
$\corr^{ab}
= 
(\corr^{ab}_{\upsilon_1 \upsilon_2})_{\upsilon_1,\upsilon_2 \in \states}$.
\end{definition}
\begin{definition}[Distance Method]
A phylogenetic reconstruction algorithm
$\acal = \{\acal_{n,k}\}_{n,k\geq 1}$ is said
to be {\em distance-based} if $\acal$ depends on the data
$(\xi^i_{[n]})_{i=1}^k \in (\states^{[n]})^k$
{\em only through the correlation matrices} $\{\corr^{ab}\}_{a,b\in [n]}$.
\end{definition}
The previous definition takes a very general view of distance-based methods:
any method that uses only pairwise sequence comparisons.
In practice, most distance-based approaches actually use a specific
{\em distance estimator}, that is, a function of $\corr^{ab}$ that converges
to $\weight(a,b)$ in probability as $n \to +\infty$.
We give two classical examples below.
\begin{example}[CFN Metric]\label{ex:cfnmetric}
In the CFN case with state space $\states=\{+,-\}$, a standard distance estimator (up to a constant) is
\begin{equation*}
\estim(\corr)
=
-\ln\left(1 - 2(\corr_{+-} + \corr_{-+})\right).
\end{equation*}
\end{example}
\begin{example}[Log-Det Distance~\cite{BarryHartigan:87, Lake:94, LoStHePe:94,Steel:94}]
More generally, a common distance estimator (up to scaling) is
the so-called {\em log-det distance}
\begin{equation*}
\estim(\corr) = -\ln|\det \corr|.
\end{equation*}
Loosely speaking, the log-det distance
can be thought as
a generalization
of the CFN metric. We will use a different generalization
of the CFN metric.
See section~\ref{section:overview}.
\end{example}

\subsection{Results}\label{sec:results}

In our main result, we prove that phylogenies under GTR models of mutation 
can be inferred using a distance-based method from $k = \poly(\log n)$.
\begin{theorem}[Main Result]\label{thm:main}
For all $\nstates \geq 2$, $0 < \quantum \leq f \leq g < g^*$ and $Q\in \rates_{\nstates}$,
there is a distance-based method 
solving the phylogenetic reconstruction problem on 
$\sphy^{f,g}_\quantum \otimes \{Q\}$
with $k = \poly(\log n)$.\footnote{
As in Theorem~\ref{thm:opt}, 
the dependence of $k$ in $\quantum$ is $\quantum^{-2}$.
}
\end{theorem}
Note that this result is a substantial improvement over Theorem~\ref{thm:essw}---at least, in a certain range of parameters---and 
that it  almost matches the bound obtained in Theorem~\ref{thm:opt}. The result is also novel in two ways over Theorem~\ref{thm:opt}:
only the distance matrix is used; the result applies to a larger class of mutation matrices.
A slightly weaker version of the result stated here appeared without proof as~\cite{Roch:08}.
Note that in~\cite{Roch:08} the result was
stated without
the discretization assumption which is in fact needed for the final step of the proof.
This is further explained in Section 7.3 of~\cite{DaMoRo:08b}.
In subsequent work~\cite{Roch:09}, 
the result stated here was improved to logarithmic sequence length, 
thereby matching Theorem~\ref{thm:opt}. 
This new result follows a similar high-level proof 
but involves stronger concentration arguments~\cite{PeresRoch:08},
as well as a simplified algorithm.

In an attempt to keep the paper as self-contained as possible 
we first give a proof in the special case of homogeneous trees.
This allows to keep the algorithmic details to a minimum. 
The proof appears in Section~\ref{section:hmg}.
We extend the result to general trees in Appendix~\ref{section:general-trees}.
The more general result relies on a combinatorial algorithm of~\cite{DaMoRo:08b}.


\subsection{Proof Overview}\label{section:overview}

\paragraph{Distance averaging.}
The basic insight behind Steel's conjecture is that
the accurate reconstruction of ancestral sequences
in the reconstruction phase can be harnessed
to perform a better reconstruction of the
phylogeny itself. For now, consider the CFN model with character space $\{+1,-1\}$ and 
assume that our phylogeny is
homogeneous with uniform branch lengths $\omega$. 
Generate $k$ i.i.d.~samples $(\s^i_V)_{i=1}^k$.
Let $a, b$ be
two internal vertices on level $h - h' < h$ (from the root). Suppose we seek to estimate the distance
between $a$ and $b$. This estimation cannot be performed directly
because the sequences at $a$ and $b$ are not known.
However, we can try to {\em estimate} these internal sequences.
Denote by $A$, $B$ the leaf
set below $a$ and $b$ respectively. 
An estimate of the sequence at $a$ is the
(properly normalized) ``site-wise average'' of the sequences
at $A$
\begin{equation}\label{eq:majority}
\bar\s_a^i = \frac{1}{|A|}\sum_{a' \in A} \frac{\s^i_{a'}}{e^{-\omega h'}},
\end{equation}
for $i=1,\ldots,k$, and similarly for $b$. 
It is not immediately clear 
how such a {\em site-wise} procedure involving {\em simultaneously}
a large number of leaves can be performed using the more aggregated information
in the correlation matrices
$\{\corr^{uv}\}_{u,v\in [n]}$.
Nevertheless, note that the quantity we are ultimately interested in
computing is the following estimate of the CFN metric between $a$ and $b$
\begin{equation*}
\bar\weight(a,b) = -\ln \left(\frac{1}{k} \sum_{i=1}^k  \bar\s_a^i\bar\s_b^i\right).
\end{equation*}
Our results are based on the following observation:
\begin{eqnarray*}
\bar\weight(a,b)
&=& -\ln \left(\frac{1}{k} \sum_{i=1}^k  \left(\frac{1}{|A|}\sum_{a' \in A} \frac{\s^i_{a'}}{e^{-\omega h'}}\right)
\left(\frac{1}{|B|}\sum_{b' \in B} \frac{\s^i_{b'}}{e^{-\omega h'}}\right)\right)\\
&=& -\ln \left(\frac{1}{|A||B|e^{-2 \omega h'}}\sum_{a' \in A}\sum_{b' \in B} 
\left(\frac{1}{k} \sum_{i=1}^k \s^i_{a'}\s^i_{b'}\right)\right)\\
&=& -\ln \left(\frac{1}{|A||B|e^{-2 \omega h'}}\sum_{a' \in A}\sum_{b' \in B} 
e^{-\hat\weight(a',b')}\right),
\end{eqnarray*}
where note that the last line depends only on distance estimates $\hat\weight(a',b')$ between leaves 
$a',b'$ in $A,B$. 
In other words, through this procedure, which we call
 {\em exponential averaging},
we perform an {\em implicit} ancestral sequence
reconstruction using only distance estimates.
One can also think of this as a variance reduction technique.
When the branch lengths are not uniform, one needs to use a
{\em weighted} version of (\ref{eq:majority}). This requires the estimation
of path lengths.

\paragraph{GTR models.}
In the case of GTR models, 
the standard log-det estimator does not lend itself well to the
exponential averaging procedure described above.
Instead, we use an estimator involving
the right eigenvector $\evr$ corresponding to the second eigenvalue $\Lambda_2$ of $Q$. 
For $a,b\in [n]$, we consider the estimator
\begin{equation}
\eweight(a,b) = -\ln \left(\evr^\top \corr^{ab} \evr\right).
\end{equation}
This choice is justified by a generalization of (\ref{eq:majority})
introduced in~\cite{MosselPeres:03}.
Note that $\evr$ may need to be estimated.

\paragraph{Concentration.}
There is a further complication in that to obtain 
results with high probability, one needs to show that
$\bar\weight(a,b)$ is {\em highly concentrated}.
However, one cannot directly apply standard concentration
inequalities because $\bar\s_a$ is {\em not bounded}.
Classical results on the reconstruction problem
imply that the variance of $\bar\s_a$ is finite---which is not quite enough.
To amplify the accuracy, we ``go down'' $\log\log n$ levels
and compute $O(\log n)$ distance estimates with
{\em conditionally independent} biases. 
By performing a majority procedure, we finally obtain a concentrated estimate.



\subsection{Organization}

In Section~\ref{section:averaging}, we provide a detailed account
of the connection between ancestral sequence reconstruction and distance
averaging.
We then give a proof of our main result in the case of homogeneous trees
in Section~\ref{section:hmg}. 
Finally, in Appendix~\ref{section:general-trees}, we give a sketch of the
proof in the general case. 

All proofs are relegated to Appendix~\ref{sec:proofs}.

%

\section{Ancestral Reconstruction and Distance Averaging}\label{section:averaging}

Let $\nstates \geq 2$, $0 < \quantum \leq f \leq g < g^* = \ln \sqrt{2}$, 
and $Q \in \rates_\nstates$ with corresponding stationary distribution $\pi > 0$.
In this section we restrict ourselves to the homogeneous case 
$\phy = \phy^{(h)} = (V,E,[n],\rt;\weight)$ where we take $h = \log_2 n$ and
$f \leq \weight_e \leq g$ and $\weight_e$ is an integer multiple of $\quantum$, 
$\forall e\in E$.
(See Examples~\ref{ex:homo} and \ref{ex:cfn} and Theorem~\ref{thm:opt}.)\footnote{Note 
that, without loss of generality, 
we can consider performing ancestral state reconstruction on a homogeneous
tree as it is always possible to ``complete'' a general tree
with zero-length edges. 
We come back to this point in Appendix~\ref{section:general-trees}.}
 
Throughout this section, we use a sequence length
$k > \log^\kappa n$ where $\kappa > 1$ is a constant 
to be determined 
later. 
We generate $k$ i.i.d.~samples $(\xi^i_{V})_{i=1}^k$ from the GTR model 
$(\phy, Q)$
with state space $\states$.

All proofs are relegated to Appendix~\ref{sec:proofs}.

\subsection{Distance Estimator}

The standard log-det estimator does not lend itself well to the
averaging procedure discussed above.
For reconstruction purposes, we instead use an estimator involving
the right eigenvector $\evr$ corresponding to the second eigenvalue $\Lambda_2$
of $Q$. For $a,b\in [n]$, consider the estimator
\begin{equation}\label{eq:eigenestim}
\eweight(a,b) = -\ln \left(\evr^\top \corr^{ab} \evr\right),
\end{equation}
where the correlation matrix $\corr^{ab}$ was introduced in Definition~\ref{def:distest}.
We first give a proof that this is indeed a legitimate distance estimator.
For more on connections between eigenvalues of the rate matrix and distance estimation, 
see e.g.~\cite{GuLi:96,GuLi:98,GrMoYa:09}.
\begin{lemma}[Distance Estimator]\label{lemma:distance}
Let $\eweight$ be as above. For all $a,b\in [n]$, we have
\begin{equation*}
\expec[e^{-\eweight(a,b)}] = e^{-\weight(a,b)}.
\end{equation*} 
\end{lemma}
For $a \in [n]$ and $i=1,\ldots,k$, let
\begin{equation*}
\sigma^i_{a} = \evr_{\xi^i_a}.
\end{equation*}
Then (\ref{eq:eigenestim}) is equivalent to
\begin{equation}\label{eq:eigenestim2}
\eweight(a,b) = -\ln \left(\frac{1}{k}\sum_{i=1}^k \sigma^i_a \sigma^i_b\right).
\end{equation}
Note that in the CFN case, we have simply $\evr = (1, -1)^\top$ and hence
(\ref{eq:eigenestim2}) can be interpreted as a generalization
of the CFN metric. 

Let
\begin{equation*}
\bar\pi = \min_\iota \pi_\iota,
\end{equation*}
and
\begin{equation*}
\bar\evr \equiv \max_i |\evr_i| \leq \frac{1}{\sqrt{\bar\pi}}.
\end{equation*}
The following lemmas show that the distance estimate above 
with sequence
length $\poly(\log n)$ is concentrated
for path lengths of order $O(\log\log n)$.
\begin{lemma}[Distorted Metric: Short Distances]\label{lem:distmet1gtr}
Let $\delta > 0$, $\eps > 0$, and $\gamma > 0$. 
There exists $\kappa > 1$, such that if the following conditions hold
for $u,v \in V$:
\begin{itemize} 
\item $\mathrm{[Small\ Diameter]}$ $\weight(u,v) < \delta\log\log(n)$,
\item $\mathrm{[Sequence\ Length]}$ $k > \log^\kappa(n)$,
\end{itemize}
then
\begin{equation*}
\left|\weight(u,v)-\eweight(u,v)\right|< \eps,
\end{equation*}
with probability at least $1-n^{-\gamma}$.
\end{lemma}
\begin{lemma}[Distorted Metric: Diameter Test]\label{lem:distmet2gtr}
Let $D > 0$, $W > 5$, $\delta > 0$, and $\gamma > 0$.
Let $u,v \in V$ not descendants of each other.
Let $u_0,v_0\in V$ be descendants
of $u,v$ respectively. 
there exists $\kappa > 1$, 
such that if the following conditions hold:
\begin{itemize} 
\item $\mathrm{[Large\ Diameter]}$ $\weight(u_0,v_0) > D + \ln W$,
\item $\mathrm{[Close\ Descendants]}$ $\weight(u_0,u), \weight(v_0,v) 
< \delta\log\log(n)$,
\item $\mathrm{[Sequence\ Length]}$ $k > \log^\kappa(n)$, 
\end{itemize}
then
\begin{equation*}
\eweight(u,v) - \weight(u_0,u) - \weight(v_0,v)
> D + \ln \frac{W}{2},
\end{equation*}
with probability at least $1-n^{-\gamma}$.
On the other hand, if the first condition above is replaced by
\begin{itemize} 
\item $\mathrm{[Small\ Diameter]}$ $\weight(u_0,v_0) 
< D + \ln \frac{W}{5}$,
\end{itemize}
then
\begin{equation*}
\eweight(u,v) - \weight(u_0,u) - \weight(v_0,v)
\leq D + \ln \frac{W}{4},
\end{equation*}
with probability at least $1-n^{-\gamma}$.
\end{lemma}

\subsection{Ancestral Sequence Reconstruction}
\label{sec:averaging}


Let $e = (x,y) \in E$
and assume that $x$ is closest to $\rt$ (in topological distance). 
We define
$\path(\rt,e) = \path(\rt,y)$,
$|e|_\rt = |\path(v,e)|$, and
\begin{equation*}
R_\rt(e) = \left(1 - \theta_e^2\right) 
\Theta_{\rt,y}^{-1},
\end{equation*}
where $\Theta_{\rt,y} = e^{-\weight(\rt,y)}$ and
$\theta_e = e^{-\weight(e)}$. 

Proposition~\ref{prop:weightedmaj} below is a variant of
Lemma 5.3 in~\cite{MosselPeres:03}. For completeness, we give a proof.
\begin{proposition}[Weighted Majority: GTR Version]\label{prop:weightedmaj}
Let $\xi_{[n]}$ be a sample from $\law_{\phy,Q}$ (see Definition~\ref{def:gtr}) 
with corresponding $\s_{[n]}$.
For a unit flow $\Psi$ from $\rt$ to $[n]$, 
consider the estimator
\begin{equation*}
S = \sum_{x \in [n]} \frac{\Psi(x) \s_x}{\Theta_{\rt,x}}.
\end{equation*}
Then, we have
\begin{equation*}
\expec[S] = 0,
\end{equation*}
\begin{equation*}
\expec[S\,|\,\s_\rt] = \s_\rt,
\end{equation*}
and
\begin{equation*}
\var[S]
= 1 + K_{\Psi},
\end{equation*}
where
\begin{equation*}
K_{\Psi} = \sum_{e \in E} R_\rt(e) \Psi(e)^2.
\end{equation*}
\end{proposition}

Let $\Psi$ be a unit flow from $\rt$ to $[n]$.
We will use the following
multiplicative decomposition of $\Psi$:
If $\Psi(x) > 0$,
we let
\begin{equation*}
\psi(e) = \frac{\Psi(y)}{\Psi(x)},
\end{equation*}
and, if instead $\Psi(x) = 0$, we let $\psi(y) = 0$.
Denoting $x_\uparrow$ the immediate ancestor of $x \in V$ and letting $\theta_x = e^{-\weight_{(x_\uparrow, x)}}$, 
it will be useful to re-write
\begin{equation}\label{eq:k}
K_{\Psi} = \sum_{h'=0}^{h-1} \sum_{x\in L^{(h)}_{h'}} (1 - \theta_x^2)\prod_{e\in \path(\rt,x)}
\frac{\psi(e)^2}{\theta_e^2},
\end{equation}
and to define the following recursion from the leaves.
For $x\in[n]$,
\begin{equation*}
K_{x,\Psi} = 0. 
\end{equation*}
Then, let $u \in V-[n]$ with children $v_1,v_2$ with corresponding edges $e_1,e_2$ and 
define
\begin{equation*}
K_{u,\Psi} = \sum_{\alpha = 1,2} ((1 - \theta_{v_\alpha}^2) + K_{v_\alpha,\Psi}) 
\left(\frac{\psi(e_\alpha)^2}{\theta_{e_\alpha}^2}\right).
\end{equation*}
Note that, from (\ref{eq:k}), we have $K_{\rt,\Psi} = K_\Psi$.
\begin{lemma}[Uniform Bound on $K_\Psi$]\label{lemma:k}
Let $\Psi$ be the homogeneous flow from $\rt$ to $[n]$. 
Then, we have 
\begin{equation*}
K_{\Psi} \leq \frac{1}{1 - e^{-2 (g^* - g)}} < +\infty.
\end{equation*}
\end{lemma}

\subsection{Distance Averaging}

The input to our tree reconstruction algorithm is the matrix
of all estimated distances between pairs of \emph{leaves} 
$\{\eweight(a,b)\}_{a,b,\in [n]}$.
For short sequences, these estimated distances are known to be accurate for leaves
that are close enough.
We now show how to compute distances between internal nodes
in a way that involves only $\{\eweight(a,b)\}_{a,b,\in [n]}$
(and previously computed internal weights)
using Proposition~\ref{prop:weightedmaj}.
However, the second-moment guarantee in Proposition~\ref{prop:weightedmaj}
is not enough to obtain good estimates with high probability.
To remedy this situation, we perform a large number 
of {\em conditionally independent} distance computations, 
as we now describe.

\begin{figure}
\begin{center}
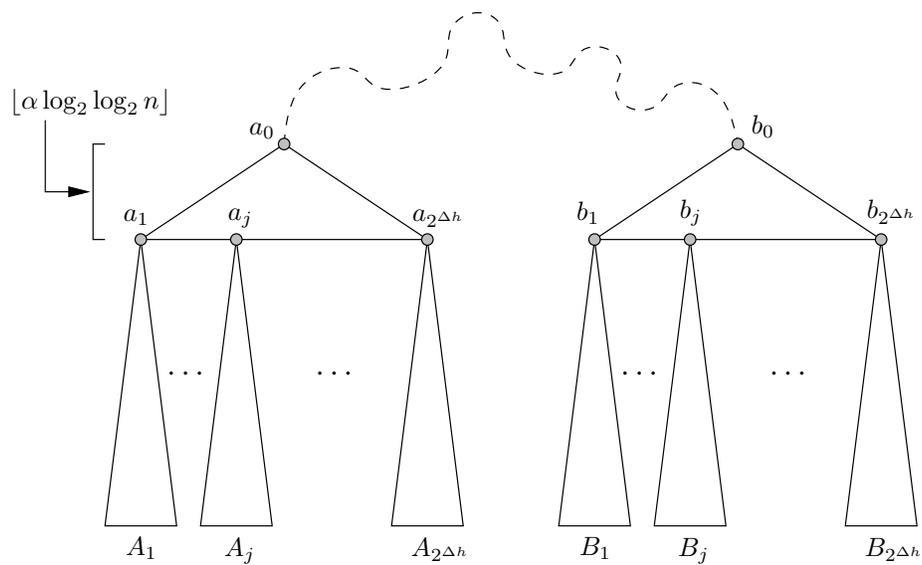\caption{Accuracy amplification.}\label{fig:amplification}
\end{center}
\end{figure}
Let $\alpha > 1$ and assume $h' > \loglog{\alpha}$.
Let $0 \leq h'' < h'$
such that
$\Delta h \equiv h' - h'' = \loglog{\alpha}$.
Let $a_0, b_0 \in L^{(h)}_{h'}$. 
For $x \in \{a,b\}$, denote by 
$x_1,\ldots,x_{2^{\Delta h}}$,
the vertices in $L^{(h)}_{h''}$ that are below $x_0$
and, for $j = 1,\ldots,2^{\Delta h}$, let $X_j$ be the leaves
of $T^{(h)}$ below $x_j$. See Figure~\ref{fig:amplification}.
Assume that we are given $\Theta_{a_0,\cdot}$, $\Theta_{b_0,\cdot}$.
For $1 \leq j \leq 2^{\Delta h}$,
we estimate $\weight(a_0,b_0)$ as follows
\begin{eqnarray*}
\bar\weight(a_j,b_j)
&\equiv& -\ln \left(\frac{1}{|A_j||B_j|}\sum_{a' \in A_j} \sum_{b' \in B_j}
\Theta^{-1}_{a_0,a'}\Theta^{-1}_{b_0,b'}
e^{-\eweight(a',b')}\right).
\end{eqnarray*}
This choice of estimator is suggested by the following observation
\begin{eqnarray*}
e^{-\bar\weight(a_j,b_j)}
&\equiv& \sum_{a' \in A_j} \sum_{b' \in B_j} 
2^{-2h''} \Theta^{-1}_{a_0,a'}\Theta^{-1}_{b_0,b'}
e^{-\eweight(a',b')}\\
&=& e^{-\weight(a_0,a_j)-\weight(b_0,b_j)} \left[\frac{1}{k} \sum_{i=1}^k \left(
\sum_{a' \in A_j} \frac{2^{-h''} \s^i_{a'}}{\Theta_{a_j,a'}}\right) 
\left( \sum_{b' \in B_j}  \frac{2^{-h''} \s^i_{b'}}{\Theta_{b_j,b'}} \right)\right].
\end{eqnarray*}
Note that the first line depends only on estimates
$(\eweight(u,v))_{u,v\in [n]}$ and $\{\Theta_{v,\cdot}\}_{v\in V_{a_0}\cup V_{b_0}}$. 
The last line is the empirical distance between
the reconstructed states at $a$ and $b$ when the flow is chosen to be uniform
in Proposition~\ref{prop:weightedmaj}. 

\begin{figure}
\begin{center}
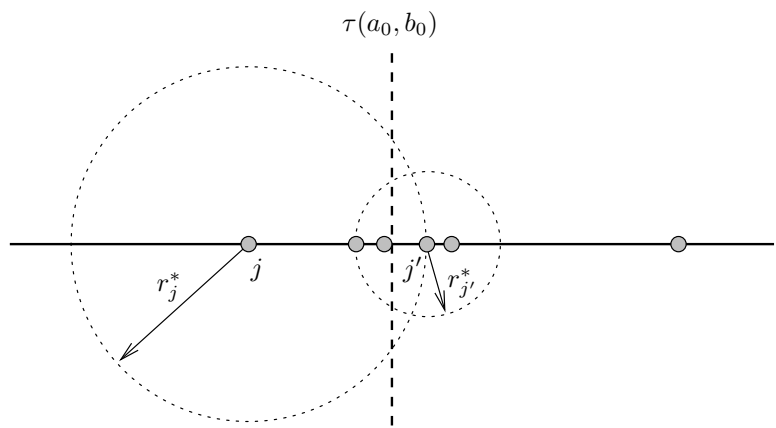\caption{Examples of dense balls with six points.
The estimate in this case would be one of the central points.}\label{fig:ball}
\end{center}
\end{figure}
For $d \in \real$ and $r > 0$, let $\ball_r(d)$ be the ball of radius $r$ around $d$.
We define the {\em dense ball around $j$} to be the smallest ball around $\bar\weight(a_{j},b_{j})$
containing at least $2/3$ of $\jcal = \{\bar\weight(a_{j'},b_{j'})\}_{j'=1}^{2^{\Delta h}}$
(as a multiset). 
The radius of the dense ball around $j$
is
\begin{equation*}
r^*_j 
= \inf\left\{r\ :\ |\ball_{r}(\bar\weight(a_{j},b_{j})) \cap \jcal| \geq \frac{2}{3} (2^{\Delta h}) \right\},
\end{equation*}
for $j = 1,\ldots, 2^{\Delta h}$.
We define our estimate of $\weight(a_0,b_0)$ to be
$\bar\weight'(a_0,b_0)
= \bar\weight(a_{j^*},b_{j^*})$,
where $j^* = \arg\min_{j} r^*_j$.
See Figure~\ref{fig:ball}.
For $D> 0$, $W > 5$, we define
\begin{equation*}
\longtest(a_0,b_0)
= \ind\left\{ 2^{-\Delta h}\left|\left\{j\,:\,
\bar\weight(a_{j},b_{j}) \leq D + \ln\frac{W}{3} \right\}\right| > \frac{1}{2}\right\}.
\end{equation*}

We extend Lemmas~\ref{lem:distmet1gtr} and~\ref{lem:distmet2gtr} to $\bar\weight'(a_0,b_0)$
and $\longtest(a_0,b_0)$.
\begin{proposition}[Deep Distance Computation: Small Diameter]\label{lem:deep1}
Let $\alpha > 1$, $D > 0$, $\gamma > 0$, and $\eps > 0$.
Let $a_0,b_0 \in L^{(h)}_{h'}$ as above. 
There exist $\kappa > 1$ 
such that if the following conditions hold:
\begin{itemize} 
\item $\mathrm{[Small\ Diameter]}$ $\weight(a_0,b_0) < D$,
\item $\mathrm{[Sequence\ Length]}$ $k > \log^\kappa(n)$,
\end{itemize}
then
\begin{equation*}
|\bar\weight'(a_0,b_0)
- \weight(a_0,b_0)| < \eps,
\end{equation*}
with probability at least $1-O(n^{-\gamma})$.
\end{proposition}

\begin{proposition}[Deep Distance Computation: Diameter Test]\label{lem:deep2}
Let $\alpha > 1$, $D > 0$, $W > 5$, and $\gamma > 0$.
Let $a_0,b_0 \in L^{(h)}_{h'}$ as above. 
There exists $\kappa > 1$ 
such that if the following conditions hold:
\begin{itemize} 
\item $\mathrm{[Large\ Diameter]}$ $\weight(a_0,b_0) > D + \ln W$,
\item $\mathrm{[Sequence\ Length]}$ $k > \log^\kappa(n)$, 
\end{itemize}
then
\begin{equation*}
\longtest(a_0,b_0) = 0,
\end{equation*}
with probability at least $1-n^{-\gamma}$.
On the other hand, if the first condition above is replaced by
\begin{itemize} 
\item $\mathrm{[Small\ Diameter]}$ $\weight(a_0,b_0) 
< D + \ln \frac{W}{5}$,
\end{itemize}
then
\begin{equation*}
\longtest(a_0,b_0) = 1,
\end{equation*}
with probability at least $1-n^{-\gamma}$.
\end{proposition}

\section{Reconstructing Homogeneous Trees}\label{section:hmg}

In this section, we prove our main result in the case of homogeneous trees.
More precisely, we prove the following.
\begin{theorem}[Main Result: Homogeneous Case]\label{thm:mainhmg}
Let $0 < \quantum \leq f \leq g < +\infty$ and denote by $\shmgphy^{f,g}_\quantum$ the set of 
all homogeneous phylogenies $\phy = (V,E,[n],\rt;\weight)$ satisfying
$f \leq \weight_e \leq g$ and $\weight_e$ is an integer multiple of $\quantum$, $\forall e\in E$.
Let $g^* = \ln \sqrt{2}$.
Then, for all $\nstates \geq 2$, $0 < \quantum \leq f \leq g < g^*$ and $Q\in \rates_{\nstates}$,
there is a distance-based method 
solving the phylogenetic reconstruction problem on 
$\shmgphy^{f,g}_\quantum \otimes \{Q\}$
with $k = \poly(\log n)$.
\end{theorem}
All proofs are relegated to Appendix~\ref{sec:proofs}.

In the homogeneous case, we can build the tree level by level using simple ``four-point''
techniques~\cite{Buneman:71}. See e.g.~\cite{SempleSteel:03,Felsenstein:04} for background
and details. See also Section~\ref{section:algorithm} below.
The underlying combinatorial algorithm we use here is essentially identical
to the one used by Mossel in~\cite{Mossel:04a}.
From Propositions~\ref{lem:deep1} and~\ref{lem:deep2},
we get that the ``local metric'' on each level is accurate as long as we compute
adequate weights. We summarize this fact in the next proposition.
For $\quantum > 0$ and $z \in \real_+$, we let $[z]_\quantum$ be the closest multiple
of $\quantum$ to $z$ (breaking ties arbitrarily).
We define
\begin{equation*}
\estdist(a_0,b_0)
=
\left\{
\begin{array}{ll}
[\bar\weight'(a_0,b_0)]_\quantum, & \text{if}\ \longtest(a_0,b_0) = 1,\\
+\infty, & \text{o.w.}
\end{array}
\right.
\end{equation*}
\begin{proposition}[Deep Distorted Metric]\label{prop:deep3}
Let $D > 0$, $W > 5$, and $\gamma > 0$.
Let $\phy = (V,E,[n],\rt;\weight) \in \shmgphy^{f,g}_\quantum$ with $g < g^*$.
Let $a_0,b_0 \in L^{(h)}_{h'}$ for $0\leq h' < h$. 
Assume we are given, for $x=a,b$, 
$\theta_{e}$ for all $e\in V_x$.
There exists $\kappa > 0$, 
such that if the following condition holds:
\begin{itemize} 
\item $\mathrm{[Sequence\ Length]}$ The sequence length is $k > \log^\kappa(n)$,
\end{itemize}
then we have, with probability at least $1 - O(n^{-\gamma})$,
\begin{equation*}
\estdist(a_0,b_0) = \weight(a_0,b_0)
\end{equation*}
under either of the following two conditions:
\begin{enumerate}
\item $\mathrm{[Small\ Diameter]}$ $\weight(a_0,b_0) < D$, or

\item $\mathrm{[Finite\ Estimate]}$ $\estdist(a_0,b_0) < +\infty$.
\end{enumerate}
\end{proposition}
It remains to show how to compute the weights, which is the purpose of the next section.

\subsection{Estimating Averaging Weights}

Proposition~\ref{prop:deep3} relies on the prior computation of the weights
$\theta_{e}$ for all $e\in V_x$, for $x=a,b$. 
In this section, we show how this estimation is performed. 

Let $a_0,b_0,c_0 \in L^{(h)}_{h'}$.
Denote by $z$ the meeting point of the paths joining $a_0, b_0, c_0$.
We define the ``three-point'' estimate
\begin{equation*}
\hat\theta_{z,a_0}
= \triplet(a_0; b_0, c_0)
\equiv
\exp\left(-\frac{1}{2}[\estdist(a_0,b_0) + \estdist(a_0,c_0) - \estdist(b_0,c_0)]\right).
\end{equation*}
Note that the expression in parenthesis is an estimate of the distance between $a_0$ and $z$.
\begin{proposition}[Averaging Weight Estimation]\label{prop:weights}
Let $a_0,b_0,c_0 \in L^{(h)}_{h'}$ as above.
Assume that the assumptions of Propositions~\ref{lem:deep1}, \ref{lem:deep2}, \ref{prop:deep3}
hold. 
Assume further that the following condition hold:
\begin{itemize} 
\item $\mathrm{[Small\ Diameter]}$ $\weight(a_0,b_0),\weight(a_0,c_0),\weight(b_0,c_0) 
< D + \ln W$,
\end{itemize}
then 
\begin{equation*}
\hat\theta_{z,a_0} = \theta_{z,a_0},
\end{equation*}
with probability at least $1-O(n^{-\gamma})$ where
$\hat\theta_{z,a_0}
= \triplet (a_0; b_0, c_0)$.
\end{proposition}

\subsection{Putting it All Together}\label{section:algorithm}

Let $0 \leq h' < h$ and
$\qcal = \{a_0,b_0,c_0,d_0\} \subseteq L^{(h)}_{h'}$.
The topology of $T^{(h)}$ restricted to $\qcal$
is completely characterized by a bi-partition or
{\em quartet split} $q$ of the form:
$a_0 b_0 | c_0 d_0$, $a_0 c_0 | b_0 d_0$ or
$a_0 d_0 | b_0 c_0$. 
The most basic operation in quartet-based reconstruction
algorithms is the inference of such quartet splits.
In distance-based methods in particular, this is usually done
by performing the so-called {\em four-point test}:
letting
\begin{equation*}
\fcal(a_0 b_0 | c_0 d_0) 
= \frac{1}{2}[\weight(a_0,c_0) + \weight(b_0,d_0) 
- \weight(a_0,b_0) - \weight(c_0,d_0)],
\end{equation*}
we have
\begin{equation*}
q 
=
\left\{
\begin{array}{ll}
a_0 b_0 | c_0 d_0 & \mathrm{if\ }\fcal(a_0,b_0|c_0,d_0) > 0\\
a_0 c_0 | b_0 d_0 & \mathrm{if\ }\fcal(a_0,b_0|c_0,d_0) < 0\\
a_0 d_0 | b_0 c_0 & \mathrm{o.w.}
\end{array}
\right.
\end{equation*} 
Of course, we cannot compute $\fcal(a_0,b_0|c_0,d_0)$
directly unless $h'=0$. Instead we use
Proposition~\ref{prop:deep3}.

\paragraph{Deep Four-Point Test.} Assume we have previously computed weights
$\theta_{e}$ for all $e\in V_{x}$,
for $x=a,b,c,d$.
We let
\begin{equation}
\overline\fcal(a_0 b_0 | c_0 d_0)
= \frac{1}{2}[\estdist(a_0,c_0) 
+ \estdist(b_0,d_0)
- \estdist(a_0,b_0) 
- \estdist(c_0,d_0)],\label{eq:fcal}
\end{equation}
and we define the {\em deep four-point test}
\begin{equation*}
\deep(a_0,b_0|c_0,d_0) = \ind\{\overline\fcal(a_0 b_0 | c_0 d_0) > f/2\},
\end{equation*} 
with $\deep(a_0,b_0|c_0,d_0) = 0$ if any of the distances
in (\ref{eq:fcal}) is infinite.
Also, we extend the {\em diameter test} $\longtest$ to
arbitrary subsets by letting $\longtest(\scal) = 1$
if and only if $\longtest(x,y) = 1$ for all
pairs $x,y \in \scal$.

\paragraph{Algorithm.}
Fix $\alpha > 1$, $D > 4g$, $W > 5$, $\gamma > 3$.
Choose $\kappa$ so as to satisfy Propositions~\ref{prop:deep3} and~\ref{prop:weights}.
Let $\zcal_{0}$ be the set of leaves.
The algorithm---a standard cherry picking algorithm---is detailed in Figure~\ref{fig:algo}.

\begin{figure*}[ht]
\framebox{
\begin{minipage}{16.6cm}
{\small \textbf{Algorithm}\\
\textit{Input:} Distance estimates $\{\eweight(a,b)\}_{a,b\in [n]}$;\\
\textit{Output:} Tree;

\begin{itemize}
\item For $h' = 1,\ldots,h-1$,
\begin{enumerate}
\item \textbf{Four-Point Test.} 
Let
\begin{equation*}
\rcal_{h'} = \{q = ab|cd\ :\ \forall a,b,c,d \in \zcal_{h'}\ \text{distinct such that}\ \deep(q) = 1\}.
\end{equation*}

\item \textbf{Cherry Picking.} Identify the cherries in $\rcal_{h'}$, that is, those pairs of vertices 
that only appear on the same side of the quartet splits in $\rcal_{h'}$. 
Let 
\begin{equation*}
\zcal_{h'+1} = \{a_1^{(h'+1)},\ldots,a_{2^{h - (h'+1)}}^{(h'+1)}\},
\end{equation*}
be the parents of the cherries in $\zcal_{h'}$

\item \textbf{Weight Estimation.} For all $z \in \zcal_{h'+1}$,
\begin{enumerate}
\item Let $x,y$ be the children of $z$.  
Choose $w$ to be any other
vertex in $\zcal_{h'}$ with $\longtest(\{x,y,w\}) = 1$.

\item Compute
\begin{equation*}
\hat\theta_{z,x} = \triplet(x;y,w).
\end{equation*}

\item Repeat the previous step interchanging the role of $x$ and $y$. 
\end{enumerate}
\end{enumerate}

\end{itemize}

}
\end{minipage}
} \caption{Algorithm.} \label{fig:algo}
\end{figure*}

\section*{Acknowledgments}

This work was triggered by a discussion with Elchanan Mossel
on lower bounds for distance methods, following a talk of
Joseph Felsenstein. In particular, Elchanan pointed out
that the distance matrix has a potentially useful correlation
structure. 

\bibliographystyle{alpha}
\bibliography{thesis}

\clearpage

\appendix

\section{Extending to General Trees}\label{section:general-trees}

It is possible to generalize the previous
arguments to general trees, using a combinatorial algorithm of~\cite{DaMoRo:08b},
thereby giving a proof of Theorem~\ref{thm:main}. 
To apply the algorithm of~\cite{DaMoRo:08b} we need 
to obtain a generalization of Proposition~\ref{prop:deep3}
for disjoint subtrees in ``general position.'' This is somewhat straightforward
and we give a quick sketch in this section.

\subsection{Basic Definitions}

The algorithm in~\cite{DaMoRo:08b} is called Blindfolded Cherry Picking.
We refer the reader to~\cite{DaMoRo:08b} for a full description
of the algorithm, which is somewhat involved.
It is very similar in spirit to the algorithm introduced in Section~\ref{section:algorithm},
except for complications due to the non-homogeneity of the tree. 
The proof in~\cite{DaMoRo:08b} is modular and relies on two main components:
a distance-based \emph{combinatorial} argument which remains unchanged in our setting; 
and a \emph{statistical} argument which we now adapt. 
The key to the latter is~\cite[Proposition 4]{DaMoRo:08b}. 
Note that~\cite[Proposition 4]{DaMoRo:08b} is \emph{not} distance-based
as it relies on a complex ancestral reconstruction function---recursive majority.
Our main contribution in this section is to show how this result
can be obtained using the techniques of the previous sections---leading
to a fully distance-based reconstruction algorithm. 

\begin{figure}
\begin{center}
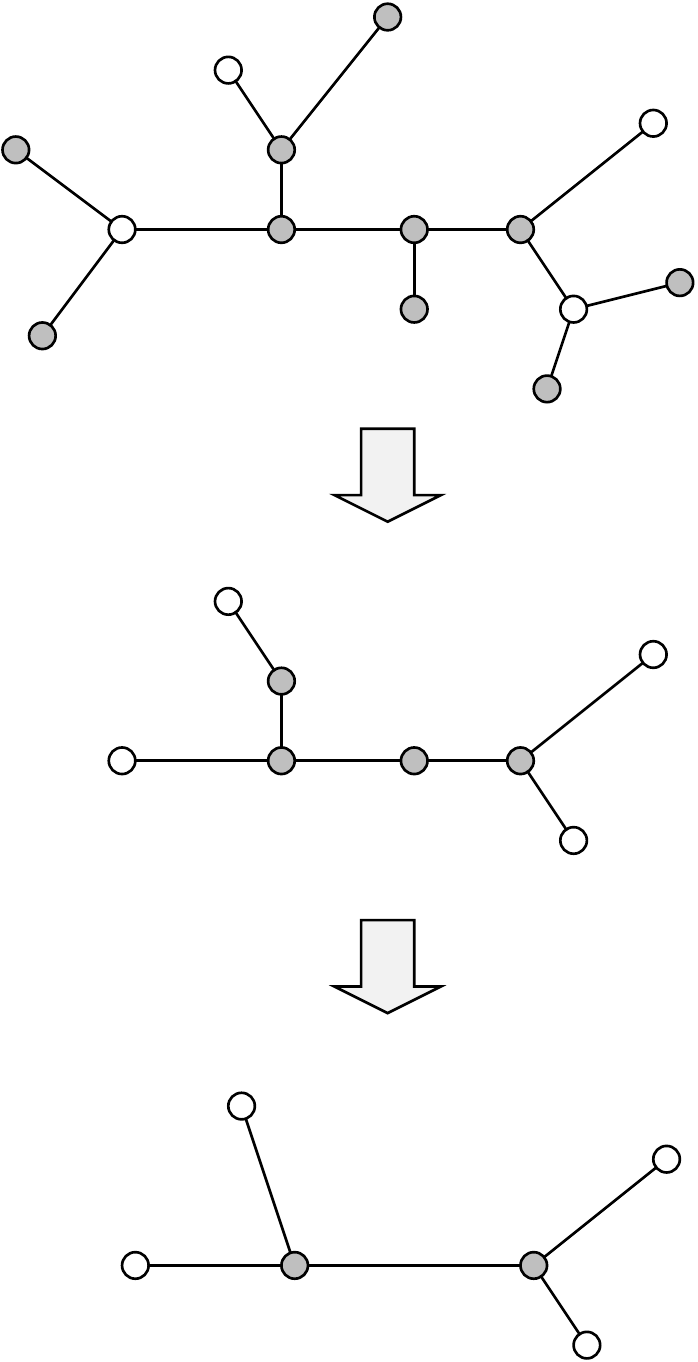\caption{Restricting the top tree to its white nodes.}\label{fig:restricted}
\end{center}
\end{figure}

\begin{figure}
\begin{center}
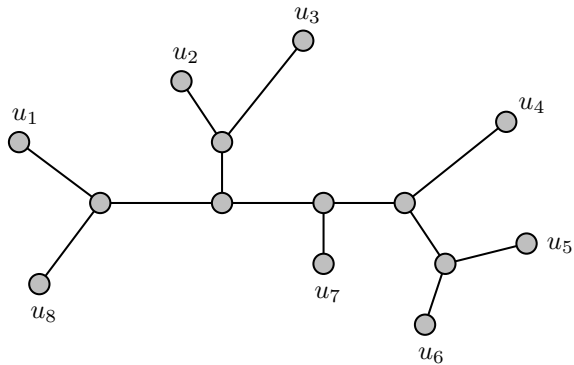\caption{
The subtrees $T|_{\{u_1,u_2,u_3,u_8\}}$
and $T|_{\{u_4,u_5,u_6,u_7\}}$ are edge-disjoint.
The subtrees $T|_{\{u_1,u_5,u_6,u_8\}}$
and $T|_{\{u_2,u_3,u_4,u_7\}}$ are edge-sharing.
}\label{fig:sharing}
\end{center}
\end{figure}

\begin{figure}
\begin{center}
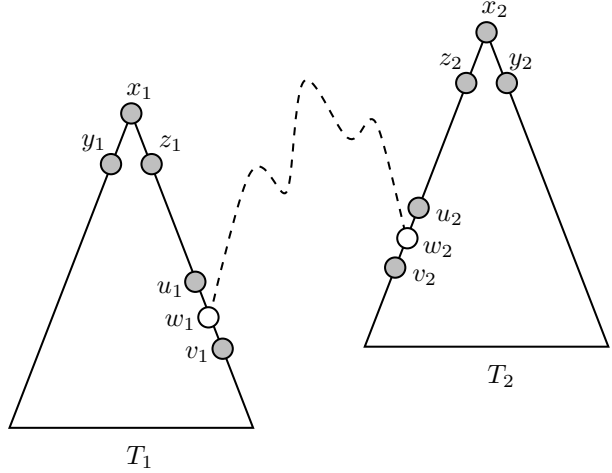\caption{
Basic Disjoint Setup (General).
The rooted subtrees $T_{1}, T_{2}$ are edge-disjoint
but are not assumed to be dangling.
The white nodes may not be in the {\em restricted} subtrees
$T_{1}, T_{2}$.
The case $w_1 = x_1$ and/or $w_2 = x_2$ is possible.
Note that if we root the tree at any node along the dashed path, the
subtrees rooted at $y_1$ and $y_2$ are edge-disjoint and dangling
(unlike $T_1$ and $T_2$).}\label{fig:remote}
\end{center}
\end{figure}

In order to explain the complications due to the non-homogeneity of the tree and
state our main result, we first need to borrow a few definitions
from~\cite{DaMoRo:08b}.
\paragraph{Basic Definitions.} 
Fix $0 < \quantum \leq f \leq g < g^*$ as in Theorem~\ref{thm:main}.
Let $\phy = (V,E,[n],\rt;\weight) \in \sphy^{f,g}_\quantum$ 
be a phylogeny with underlying tree $T = (V,E)$.
In this section, we sometimes refer to the edge set, vertex set and 
leaf set of a tree $T'$ as $\ecal(T')$, $\vcal(T')$, and $\lcal(T')$ respectively.
\begin{definition}[Restricted Subtree]
Let $V' \subseteq V$ be a subset of the vertices of $T$.
The {\em subtree of $T$ restricted to $V'$} is the tree $T'$
obtained by 1) keeping only nodes and edges on paths between
vertices in $V'$ and 2) by then contracting all paths composed
of vertices of degree 2, except the nodes in $V'$. We sometimes use the notation
$T' = T|_{V'}$. See Figure~\ref{fig:restricted} for an example.
\end{definition}
\begin{definition}[Edge Disjointness] \label{def:disjoint}
Denote by $\path_T(x,y)$ the path (sequence of edges) connecting
$x$ to $y$ in $T$. We say that two restricted subtrees $T_1, T_2$ of $T$ are
{\em edge disjoint} if
\[
\path_T(x_1,y_1) \cap \path_T(x_2,y_2) = \emptyset,
\]
for all $x_1, y_1 \in \lcal(T_1)$ and $x_2, y_2 \in \lcal(T_2)$. We
say that $T_1,T_2$ are {\em edge sharing} if they are not edge
disjoint. See Figure~\ref{fig:sharing} for an example.
\end{definition}
\begin{definition}[Legal Subforest]
We say that a tree is a rooted full binary tree if all its internal
nodes have degree 3 except the root which has degree 2. A restricted
subtree $T_1$ of $T$ is a {\em legal subtree} of $T$ if it is also
a rooted full binary tree. We say that a forest 
\begin{equation*}
\forestno = \{T_1,T_2,\ldots \},
\end{equation*}
is {\em legal subforest} of $T$ if the $T_\onetwo$'s are
\emph{edge-disjoint} legal subtrees of $T$.
We denote by $\rho(\forestno)$ the set of roots of $\forestno$.
\end{definition}
\begin{definition}[Dangling Subtrees]
We say that two edge-disjoint legal subtrees $T_1$, $T_2$ of $T$
are {\em dangling} if there is a choice of root for $T$
\emph{not in $T_1$ or $T_2$}
that is consistent with the rooting of both $T_1$ and $T_2$.
See Figure~\ref{fig:remote} below for an example where two legal, edge-disjoint subtrees
are \emph{not} dangling.
\end{definition}
\begin{definition}[Basic Disjoint Setup (General)]\label{def:bds}
Let $T_1 = T_{x_1}$ and $T_2 = T_{x_2}$ be two restricted subtrees of $T$ rooted
at $x_1$ and $x_2$ respectively.
Assume further that
$T_1$ and $T_2$ are {\em edge-disjoint}, but not necessarily {\em dangling}.
Denote by $y_\onetwo, z_\onetwo$ the children of $x_\onetwo$
in $T_\onetwo$, 
$\onetwo=1,2$. 
Let $w_\onetwo$ be the node in $T$ where the path
between $T_1$ and $T_2$ meets $T_\onetwo$, $\onetwo = 1,2$.
Note that $w_\onetwo$ may not be in $T_\onetwo$ since $T_\onetwo$ is {\em restricted}, $\onetwo = 1,2$. 
If $w_\onetwo \neq x_\onetwo$, assume without loss of generality that $w_\onetwo$
is in the subtree of $T$ rooted at $z_\onetwo$, $\onetwo = 1,2$.
We call this configuration the \emph{Basic Disjoint Setup
(General)}.
See
Figure~\ref{fig:remote}. 
Let $\weight(T_1,T_2)$ be the length of the path between $w_1$ and
$w_2$ in the metric $\weight$.
\end{definition}

\subsection{Deep Distorted Metric} 

Our reconstruction algorithm for homogeneous
trees (see Section~\ref{section:hmg}) builds the tree level by level and only encounters situations
where one has to compute the distance between two \emph{dangling} subtrees
(that is, the path connecting the subtrees ``goes above them''). However,
when reconstructing general trees by growing a subforest from the leaves, 
more general situations such as the one depicted in 
Figure~\ref{fig:remote} cannot be avoided and have to be dealt with carefully.

Hence, our goal in this subsection 
is to compute the distance between the internal nodes $x_1$ and
$x_2$ in the Basic Disjoint Setup (General). 
We have already shown
how to perform this computation when $T_1$ and $T_2$ are {\em dangling}, as
this case is handled easily by Proposition~\ref{prop:deep3} (after a slight modification
of the distance estimate; see below).
However, in the general case depicted in Figure~\ref{fig:remote}, there is a complication.
When $T_1$ and $T_2$ are {\em not} dangling, the reconstructed sequences
at $x_1$ and $x_2$
are {\em not} conditionally independent.
But it can be shown that for the algorithm Blindfolded Cherry Picking
to work properly, we need: 1) to compute the distance between $x_1$
and $x_2$ correctly when the two subtrees are close and dangling;
2) detect when the two subtrees are far apart (but an accurate distance
estimate is not required in that case). This turns out to be enough because
the algorithm
Blindfolded Cherry Picking ensures roughly that close reconstructed subtrees 
are always dangling.
We refer the reader to~\cite{DaMoRo:08b} for details.  

The key point is the following: if one computes the distance
between $y_1$ and $y_2$ {\em rather than} the distance between $x_1$ and $x_2$, 
then the dangling assumption
is satisfied (re-root the tree at any node 
along the path connecting $w_1$ and $w_2$).
However, when the algorithm has only reconstructed $T_1$ and $T_2$, 
we cannot tell which pair in $\{y_1, z_1\}\times\{y_2, z_2\}$ 
is the right one to use for the distance estimation.
Instead, we compute the distance for all pairs in
$\{y_1, z_1\}\times\{y_2, z_2\}$
and the following then holds: in the dangling case, all these distances will
agree (after subtracting the length of the edges between $x_1, x_2$ and
$\{y_1, z_1, y_2, z_2\}$); in the general case, at least one is correct. 
This is the basic observation behind the routine \distmet\ in Figure~\ref{fig:distmet}
and the proof of Proposition~\ref{prop:deepgeneral} below.
We slightly modify the definitions of Section~\ref{section:hmg}. 

Using the notation of Definition~\ref{def:bds},
fix $(a_0,b_0) \in \{y_1,z_1\}\times\{y_2,z_2\}$.
For $v\in V$ and $\ell\in L$, 
denote by 
$|\ell|_v$ be the graph distance (that is, the number of edges)
between $v$ and leaf $\ell$.
Assume that we are given $\theta_e$ for all $e\in \ecal(T_{a_0})\cup\ecal(T_{b_0})$.
Let $\alpha > 1$ and $\Delta h = \loglog{\alpha}$.
Imagine (minimally) completing the subtrees below $a_0$ and $b_0$ with
$0$-length edges so that the leaves below $a_0$ and $b_0$ are at distance
at least $\Delta h$.
For $x \in \{a,b\}$, denote by 
$x_1,\ldots,x_{2^{\Delta h}}$,
the vertices below $x_0$ at distance $\Delta h$ from $x_0$
and, for $j = 1,\ldots,2^{\Delta h}$, let $X_j$ be the leaves
of $T$ below $x_j$.
For $1 \leq j \leq 2^{\Delta h}$,
we estimate $\weight(a_0,b_0)$ as follows
\begin{eqnarray*}
\bar\weight(a_j,b_j)
&\equiv& -\ln \left(\frac{1}{|A_j||B_j|}\sum_{a' \in A_j} \sum_{b' \in B_j}
\Theta^{-1}_{a_0,a'}\Theta^{-1}_{b_0,b'}
e^{-\eweight(a',b')}\right).
\end{eqnarray*}
Note that, because the tree is binary, it holds that
\begin{equation*}
\sum_{a' \in A_j} \sum_{b' \in B_j} 2^{-|a'|_{a_j} - |b'|_{b_j}} 
= \sum_{a' \in A_j} 2^{-|a'|_{a_j}} \sum_{b' \in B_j}  2^{- |b'|_{b_j}} = 1,
\end{equation*}
and we can think of the weights on $A_j$ (similarly for $B_j$) as resulting from a homogeneous flow 
$\Psi_{a_j}$ from $a_j$ to $A_j$. Then, the bound on the variance of
\begin{equation*}
S_{a_j} \equiv \sum_{a' \in A_j} 2^{-|a'|_{a_j}}
\Theta^{-1}_{a_j,a'}\s_{a'},
\end{equation*}
in Proposition~\ref{prop:weightedmaj}
still holds with
\begin{equation*}
K_{a_j,\Psi_{a_j}} = \sum_{e\in \ecal(T_{a_j})} R_{a_j}(e) \Psi(e)^2.
\end{equation*}
Moreover $K_{a_j,\Psi_{a_j}}$ is uniformly bounded following
an argument identical to (\ref{eq:kbound}) in the proof of 
Lemma~\ref{lemma:k}. 
For $d \in \real$ and $r > 0$, let $\ball_r(d)$ be the ball of radius $r$ around $d$.
We define the {\em dense ball around $j$} to be the smallest ball around $\bar\weight(a_{j},b_{j})$
containing at least $2/3$ of $\jcal = \{\bar\weight(a_{j'},b_{j'})\}_{j'=1}^{2^{\Delta h}}$
(as a multiset). 
The radius of the dense ball around $j$
is
\begin{equation*}
r^*_j 
= \inf\left\{r\ :\ |\ball_{r}(\bar\weight(a_{j},b_{j})) \cap \jcal| \geq \frac{2}{3} (2^{\Delta h}) \right\},
\end{equation*}
for $j = 1,\ldots, 2^{\Delta h}$.
We define our estimate of $\weight(a_0,b_0)$ to be
$\bar\weight'(a_0,b_0)
= \bar\weight(a_{j^*},b_{j^*})$,
where $j^* = \arg\min_{j} r^*_j$.
See Figure~\ref{fig:ball}.
For $D> 0$, $W > 5$, we define
\begin{equation*}
\longtest(a_0,b_0)
= \ind\left\{ 2^{-\Delta h}\left|\left\{j\,:\,
\bar\weight(a_{j},b_{j}) \leq D + \ln\frac{W}{3} \right\}\right| > \frac{1}{2}\right\},
\end{equation*}
and we let
\begin{equation*}
\estdist(a_0,b_0)
=
\left\{
\begin{array}{ll}
[\bar\weight'(a_0,b_0)]_\quantum, & \text{if}\ \longtest(a_0,b_0) = 1,\\
+\infty, & \text{o.w.}
\end{array}
\right.
\end{equation*}

\begin{figure*}[ht]
\framebox{
\begin{minipage}{16.6cm}
{\small \textbf{Algorithm} \distmet\\
\textit{Input:} 
Rooted forest $\forestno = \{T_1,T_2\}$ rooted at vertices $x_1, x_2$; 
weights $\weight_e$, for all $e\in \ecal(T_1)\cup\ecal(T_2)$;\\
\textit{Output:} Distance $\Upsilon$;
\begin{itemize}
\item \itemname{Children} Let
$y_\onetwo, z_\onetwo$ be the children of $x_\onetwo$ in $\fcal$
for $\onetwo = 1,2$ (if $x_\onetwo$ is a leaf,
set $z_\onetwo = y_\onetwo = x_\onetwo$);

\item $\mathrm{[Distance\ Computations]}$ For all pairs 
$(a,b) \in \{y_1, z_1\}\times\{y_2, z_2\}$, compute 
\begin{equation*}
\metricno(a,b) 
:= \estdist(a,b) - \weight(a,x_1) - \weight(b,x_2);
\end{equation*}

\item \itemname{Multiple\ Test} If
\begin{eqnarray*}
&&\max\{\left|\metricno(r_1^{(1)}, r_2^{(1)}) - \metricno(r_1^{(2)}, r_2^{(2)})\right|\ :\\
&&\qquad\qquad (r_1^{(\onetwo)}, r_2^{(\onetwo)}) \in \{y_1, z_1\}\times\{y_2, z_2\}, \onetwo = 1,2\} 
= 0,
\end{eqnarray*}
return $\Upsilon := \metricno(z_1, z_2)$,
otherwise return $\Upsilon := +\infty$
(return $\Upsilon := +\infty$ if any of the distances above is $+\infty$).
\end{itemize}
}
\end{minipage}
} \caption{Routine \distmet.} \label{fig:distmet}
\end{figure*}
\begin{proposition}[Accuracy of \distmet]\label{prop:deepgeneral}
Let $\alpha > 1$, $D > 0$, $W > 5$, $\gamma > 0$ and $g < g' < g^*$.
Consider the Basic Disjoint Setup (General) with $\forestno = \{T_1, T_2\}$
and $\quartet = \{y_1,z_1,y_2,z_2\}$.
Assume we are given  
$\theta_{e}$ for all $e\in \ecal(T_1)\cup\ecal(T_2)$.
Let $\Upsilon$ denote the output of \distmet\ in Figure~\ref{fig:distmet}.
There exists $\kappa > 0$, 
such that if the following condition holds:
\begin{itemize} 
\item $\mathrm{[Edge\ Length]}$ It holds that $\weight(e)\leq g'$, 
$\forall e \in \ecal(T_{x})$, $x\in\quartet$\footnote{For technical reasons explained
in~\cite{DaMoRo:08b}, we allow edges slightly longer than the upper bound $g$.};

\item $\mathrm{[Sequence\ Length]}$ The sequence length is $k > \log^\kappa(n)$,

\end{itemize}
then we have, with probability at least $1 - O(n^{-\gamma})$,
\begin{equation*}
\Upsilon = \weight(x_1,x_2)
\end{equation*}
under either of the following two conditions:
\begin{enumerate}
\item $\mathrm{[Dangling\ Case]}$ $T_1$ and $T_2$ are dangling
and $\weight(T_1,T_2) < D$, or

\item $\mathrm{[Finite\ Estimate]}$ $\Upsilon < +\infty$.
\end{enumerate}
\end{proposition}
The proof, which is a simple combination of the proof of Proposition~\ref{prop:deep3} and
the remarks above the statement of Proposition~\ref{prop:deepgeneral}, is left out.

\paragraph{Full Algorithm.}
The rest of the Blindfolded Cherry Picking algorithm is unchanged except
for an additional step to compute averaging weights as in the algorithm
of Section~\ref{section:hmg}. This concludes our sketch of the proof
of Theorem~\ref{thm:main}.

\section{Proofs}\label{sec:proofs}

\begin{app-proof}{Lemma~\ref{lemma:distance}}
Note that $\expec[\corr^{ab}_{ij}] = \pi_i \left(e^{-\weight(a,b) Q}\right)_{ij}$.
Then
\begin{eqnarray*}
\expec\left[\evr^\top\corr^{ab} \evr\right] 
&=& \sum_{i\in\states} \evr_i \sum_{j\in\states}\pi_i \left(e^{-\weight(a,b) Q}\right)_{ij} \evr_j\\
&=& \sum_{i\in\states} \evr_i (\pi_i e^{-\weight(a,b)}\evr_i)\\
&=& e^{-\weight(a,b)} \sum_{i\in\states} \pi_i \evr_i^2\\
&=& e^{-\weight(a,b)}.
\end{eqnarray*}
\end{app-proof}

\begin{app-proof}{Lemma~\ref{lem:distmet1gtr}}
By Azuma's inequality we get
\begin{eqnarray*}
&&\prob\left[\eweight(u,v) > \weight(u,v) + \eps\right]\\
&& \qquad = \prob\left[\frac{1}{k} \sum_{i=1}^k \s_u^i \s_v^i < e^{- \weight(u,v) - \eps}\right]\\
&& \qquad = \prob\left[\frac{1}{k} \sum_{i=1}^k \s_u^i \s_v^i < e^{- \weight(u,v)} - (1 - e^{- \eps})e^{- \weight(u,v)}\right]\\
&& \qquad \leq \prob\left[\frac{1}{k} \sum_{i=1}^k \s_u^i \s_v^i < 
\expec[\s_u^1 \s_v^1] - (1 - e^{- \eps})e^{- \delta\log\log(n)}\right]\\
&& \qquad \leq \exp\left(-\frac{\left((1 - e^{- \eps})e^{- \delta\log\log(n)}\right)^2}
{2k(2\bar\evr/k)^2}\right)\\
&& \qquad \leq n^{-\gamma},
\end{eqnarray*}
for $\kappa$ large enough depending on
$\delta,\eps,\gamma$. Above, we used that 
\begin{equation*}
\weight(u,v) = - \ln \expec[\s_u^1 \s_v^1].
\end{equation*}
A similar inequality holds for the other direction.
\end{app-proof}

\begin{app-proof}{Lemma~\ref{lem:distmet2gtr}}
Assume the first three conditions hold.
By Azuma's inequality we get
\begin{eqnarray*}
&&\prob\left[\eweight(u,v) - \weight(u_0,u) - \weight(v_0,v)
\leq D + \ln \frac{W}{2}\right]\\
&& \qquad =\prob\left[\frac{1}{k} \sum_{i=1}^k \s_u^i \s_v^i 
\geq 2W^{-1} e^{- D}e^{- \weight(u_0,u) - \weight(v_0,v)}\right]\\
&& \qquad\leq\prob\left[\frac{1}{k} \sum_{i=1}^k \s_u^i \s_v^i 
\geq e^{- \weight(u,v)} + W^{-1} e^{- D}e^{- \weight(u_0,u) - \weight(v_0,v)}\right]\\
&& \qquad = \prob\left[\frac{1}{k} \sum_{i=1}^k \s_u^i \s_v^i \geq 
\expec[\s_u^1 \s_v^1] + W^{-1} e^{- D}e^{- \weight(u_0,u) - \weight(v_0,v)}\right]\\
&& \qquad\leq \exp\left(-\frac{\left(W^{-1} e^{- D}e^{- \weight(u_0,u) - \weight(v_0,v)}\right)^2}
{2k(2\bar\evr/k)^2}\right)\\
&& \qquad\leq  n^{-\gamma},
\end{eqnarray*}
for $\kappa$ large enough.
A similar argument gives the second claim.
\end{app-proof}

\begin{app-proof}{Proposition~\ref{prop:weightedmaj}}
We follow the proofs of~\cite{EvKePeSc:00, MosselPeres:03}.
Let $\unit_i$ be the unit vector in direction $i$.
Let $x\in [n]$, then
\begin{equation*}
\expec[\unit_{\xi_x}^\top\,|\,\xi_\rt] 
= \unit^\top_{\xi_\rt} e^{\weight(\rt,x)Q}.
\end{equation*}
Therefore,
\begin{equation*}
\expec[\s_x\,|\,\s_\rt] 
= \unit^\top_{\xi_\rt} e^{\weight(\rt,x)Q} \evr
= \s_\rt e^{-\weight(\rt,x)},
\end{equation*}
and
\begin{equation*}
\expec[S\,|\,\s_\rt] 
= \sum_{x \in [n]} \frac{\Psi(x) \s_\rt e^{-\weight(\rt,x)}}{\Theta_{\rt,x}}
= \s_\rt \sum_{x \in [n]} \Psi(x)
= \s_\rt.
\end{equation*}
In particular,
\begin{equation*}
\expec[S] 
= \sum_{\iota\in\states} \pi_i \evr_i = 0.
\end{equation*}

For $x,y \in [n]$, let $x\land y$ be the meeting point of the paths
between $\rt,x,y$. We have
\begin{eqnarray*}
\expec[\s_x\s_y]
&=& \sum_{\iota\in\states} \prob[\xi_{x\land y} = \iota] \expec[\s_x\s_y\,|\,\xi_{x\land y} = \iota]\\
&=& \sum_{\iota \in \states} \pi_{\iota} \expec[\s_x\,|\,\xi_{x\land y} = \iota]
\expec[\s_y\,|\,\xi_{x\land y} = \iota]\\
&=& \sum_{\iota \in \states} \pi_{\iota} e^{-\weight(x\land y, x)} \evr_{\iota}
e^{-\weight(x\land y, y)} \evr_{\iota}\\
&=& e^{-\weight(x,y)} \sum_{\iota \in \states} \pi_{\iota}  \evr_{\iota}^2\\
&=& e^{-\weight(x,y)}.
\end{eqnarray*}
Then
\begin{eqnarray*}
\var[S] &=& \expec[S^2]\\
&=& \sum_{x,y\in [n]} \frac{\Psi(x)\Psi(y)}{\Theta_{\rt,x}\Theta_{\rt,y}} \expec[\s_x\s_y]\\
&=& \sum_{x,y\in [n]} \Psi(x)\Psi(y)e^{2\weight(\rt, x\land y)}.
\end{eqnarray*}
For $e \in E$, let $e = (e_\uparrow,e_{\downarrow})$ where $e_\uparrow$ is the vertex
closest to $\rt$. Then, by a telescoping sum, for $u \in V$
\begin{eqnarray*}
\sum_{e \in \path(\rt,u)} R_{\rt}(e)
&=& \sum_{e \in \path(\rt,u)} e^{2\weight(\rt, e_\downarrow)}
- \sum_{e \in \path(\rt,u)} e^{2\weight(\rt, e_\uparrow)}\\
&=& e^{2\weight(\rt,u)} - 1,
\end{eqnarray*}
and therefore
\begin{eqnarray*}
\expec[S^2]
&=& \sum_{x,y\in [n]} \Psi(x)\Psi(y)e^{2\weight(v, x\land y)}\\
&=& \sum_{x,y\in [n]} \Psi(x)\Psi(y) \left(1 + \sum_{e \in \path(\rt,x \land y)} R_{\rt}(e)\right)\\
&=& 1 + \sum_{e\in E} R_\rt(e) \sum_{x,y\in [n]} \ind\{e\in \path(\rt,x\land y)\}  \Psi(x)\Psi(y)\\
&=& 1 + \sum_{e\in E} R_\rt(e) \Psi(e)^2.
\end{eqnarray*}
\end{app-proof}

\begin{app-proof}{Lemma~\ref{lemma:k}}
From (\ref{eq:k}), we have
\begin{eqnarray}
K_{\rt,\Psi}
&\leq& \sum_{i=0}^{h-1} (1 - e^{-2g}) 2^{h-i} \frac{e^{2 (h-i) g}}{2^{2 (h-i)}}\nonumber\\
&\leq& \sum_{j=1}^{h} e^{2 j g} e^{-(2 \ln \sqrt{2})j}\nonumber\\
&=& \sum_{j=1}^{h} e^{2 j (g - g^*)}\nonumber\\
&\leq& \sum_{j=0}^{+\infty} (e^{-2 (g^* - g)})^j\nonumber\\
&=& \frac{1}{1 - e^{-2 (g^* - g)}} < +\infty,\label{eq:kbound}
\end{eqnarray}
where recall that $g^* = \ln\sqrt{2}$.
\end{app-proof}

\begin{app-proof}{Proposition~\ref{lem:deep1}}
Let 
\begin{equation*}
\zcal = \{a_j\}_{j=1}^{2^{\Delta h}}
\cup \{b_j\}_{j=1}^{2^{\Delta h}},
\end{equation*}
and
\begin{equation*}
\ecal = \{(\s^i_{\zcal})_{i=1}^k\}.
\end{equation*}
For $j = 1, \ldots, 2^{\Delta h}$, 
let
\begin{equation*}
e^{-\eweight(a_j,b_j)} 
= \frac{1}{k}\sum_{i=1}^k \s^i_{a_j} \s^i_{b_j}.
\end{equation*}
Note that
\begin{equation*}
\expec[e^{-\eweight(a_j,b_j)}] = e^{-\weight(a_j,b_j)},
\end{equation*}
by Lemma~\ref{lemma:distance}.
For $i = 1, \ldots, k$, $j = 1, \ldots, 2^{\Delta h}$, and $x \in \{a,b\}$, 
let
\begin{equation*}
\bar\s^i_{x_j} = \sum_{x' \in X_j} \frac{2^{-h''} \s^i_{x'}}{\Theta_{x_j, x'}}.
\end{equation*}
By the Markov property, it follows that, conditioned on $\ecal$,
\begin{equation*}
\{\bar\s^i_{x_j}\ :\ i = 1, \ldots, k,\ 
j = 1, \ldots, 2^{\Delta h},\ 
x \in \{a,b\}\},
\end{equation*}
are mutually independent. Moreover, by Proposition~\ref{prop:weightedmaj}, we have
\begin{equation*}
\expec[\bar\s^i_{x_j}\,|\,\ecal] 
= \s^i_{x_j},
\end{equation*}
and
\begin{equation*}
\expec[(\bar\s^i_{x_j})^2\,|\,\ecal] \leq 2\bar\pi^{-1} \frac{1}{1 - e^{-2(g^* - g)}}.
\end{equation*}
Therefore, for any $\zeta > 0$ there exists $\kappa > 1$ such that
\begin{eqnarray*}
\expec[e^{-\bar\weight(a_j,b_j)}\,|\,\ecal]
&=& e^{\weight(a_0,a_j) + \weight(b_0,b_j)}\left(\frac{1}{k} \sum_{i=1}^k \expec[\bar\s^i_{a_j} \bar\s^i_{b_j}\,|\,\ecal]\right)\\
&=& e^{\weight(a_0,a_j) + \weight(b_0,b_j)}\left(\frac{1}{k} \sum_{i=1}^k \expec[\bar\s^i_{a_j}\,|\,\ecal]
\expec[\bar\s^i_{b_j}\,|\,\ecal]\right)\\
&=& e^{-(\eweight(a_j,b_j) - \weight(a_j,a_0) - \weight(b_j,b_0))},
\end{eqnarray*}
and
\begin{eqnarray*}
\var[e^{-\bar\weight(a_j,b_j)}\,|\,\ecal]
&=& e^{2\weight(a_0,a_j) + 2\weight(b_0,b_j)}\left(\frac{1}{k^2} \sum_{i=1}^k \var[\bar\s^i_{a_j} \bar\s^i_{b_j}\,|\,\ecal]\right)\\
&\leq& e^{2\weight(a_0,a_j) + 2\weight(b_0,b_j)}\left(\frac{1}{k^2} \sum_{i=1}^k \expec[(\bar\s^i_{a_j})^2\,|\,\ecal]
\expec[(\bar\s^i_{b_j})^2\,|\,\ecal]\right)\\
&\leq& \frac{e^{2\weight(a_0,a_j) + 2\weight(b_0,b_j)}}{k}\left(\frac{2}{1 - e^{-2(g^* - g)}}\right)^2\\
&\leq& \frac{e^{4 g\loglog{ \alpha}}}{k}\left(\frac{2}{1 - e^{-2(g^* - g)}}\right)^2\\
&\leq& \zeta.
\end{eqnarray*}

Take $\kappa$ large enough such that Lemma~\ref{lem:distmet1gtr} holds for diameter
$2 g\loglog{ \alpha} + D$, precision $\eps/6$
and failure probability $O(n^{-(\gamma+1)})$. 
Then, we have
\begin{equation}\label{eq:eweight}
|\eweight(a_j,b_j) - \weight(a_j,b_j)| < \eps/6,
\quad \forall j=1,\ldots,2^{\Delta h},
\end{equation}
with probability $1 - O(n^{-\gamma})$.
Let $\ecal'$ be the event that (\ref{eq:eweight}) holds.
Note that
\begin{equation*}
\expec[e^{-\bar\weight(a_j,b_j)}\,|\,\ecal\cap \ecal']
= e^{-(\eweight(a_j,b_j) - \weight(a_j,a_0) - \weight(b_j,b_0))}.
\end{equation*}
Let
\begin{equation*}
\eps' = \min\{(e^{\eps/3} - e^{\eps/6})e^{-D}, (e^{-\eps/6} - e^{-\eps/3})e^{-D}\}.
\end{equation*}
By Chebyshev's inequality, we have that
\begin{eqnarray*}
&&\prob\left[
\bar\weight(a_j,b_j)
< \weight(a_0,b_0) 
- \frac{1}{3}\eps
\,|\,\ecal\cap\ecal'\right]\\
&& \qquad \leq
\prob\left[
e^{-\bar\weight(a_j,b_j)}
>
e^{-\weight(a_0,b_0) + \eps/3}
\,|\,\ecal\cap\ecal'\right]\\
&& \qquad \leq
\prob\left[
e^{-\bar\weight(a_j,b_j)}
- \expec[e^{-\bar\weight(a_j,b_j)}\,|\,\ecal\cap \ecal']
> [e^{-\weight(a_0,b_0)+\eps/3} - e^{-(\eweight(a_j,b_j) - \weight(a_j,a_0) - \weight(b_j,b_0))}]
\,|\,\ecal\cap\ecal'\right]\\
&& \qquad \leq
\prob\left[
e^{-\bar\weight(a_j,b_j)}
- \expec[e^{-\bar\weight(a_j,b_j)}\,|\,\ecal\cap \ecal']
> e^{-\weight(a_0,b_0)}[e^{\eps/3} - e^{\weight(a_j,b_j)-\eweight(a_j,b_j)}]
\,|\,\ecal\cap\ecal'\right]\\
&& \qquad \leq
\prob\left[
e^{-\bar\weight(a_j,b_j)}
- \expec[e^{-\bar\weight(a_j,b_j)}\,|\,\ecal\cap \ecal']
> e^{-D}[e^{\eps/3} - e^{\eps/6}]
\,|\,\ecal\cap\ecal'\right]\\
&& \qquad \leq
\prob\left[
e^{-\bar\weight(a_j,b_j)}
- \expec[e^{-\bar\weight(a_j,b_j)}\,|\,\ecal\cap \ecal']
> \eps'
\,|\,\ecal\cap\ecal'\right]\\
&& \qquad \leq 1/12, 
\end{eqnarray*}
for $\zeta$ small enough.
A similar argument holds for
\begin{equation*}
\prob\left[
\bar\weight(a_j,b_j)
> \weight(a_0,b_0) 
+ \frac{1}{3}\eps
\,|\,\ecal\cap\ecal'\right]
\leq 1/12.
\end{equation*}
Conditioned on $\ecal\cap\ecal'$,
\begin{equation*}
\lcal = 
2^{-\Delta h}\sum_{j=1}^{2^{\Delta h}}\ind\left\{|\bar\weight(a_{j},b_{j}) 
- \weight(a_{0},b_{0})| < \frac{1}{3}\eps\right\}
\end{equation*}
is a sum of independent $\{0,1\}$-variables with average at least $5/6$.
By Azuma's inequality, we have
\begin{eqnarray*}
\prob[\lcal  
\leq 2/3 \,|\, \ecal\cap\ecal']
&\leq& \prob[\lcal - \expec[\lcal\,|\,\ecal\cap\ecal'] 
< - 1/6 \,|\, \ecal\cap\ecal']\\
&\leq& \exp\left(- \frac{(1/6)^2}{2 (2^{-\Delta h})^2 2^{\Delta h}}\right)\\
&\leq& \exp\left(-O(\LogLog{\alpha}) \right)\\
&\leq& O(n^{-\gamma}),
\end{eqnarray*}
where we used that $\alpha > 1$.

This implies that
with (unconditional) probability at least $1 - O(n^{-\gamma})$
\begin{equation*}
\left|\left\{|\bar\weight(a_{j},b_{j}) 
- \weight(a_{0},b_{0})| < \frac{1}{3}\eps\right\}\right|
> \frac{2}{3} (2^{\Delta h}).
\end{equation*}
In particular, there exists $j$ such that
\begin{equation*}
\left|\left\{j'\,:\,|\bar\weight(a_{j'},b_{j'}) 
- \bar\weight(a_{j},b_{j})| < \frac{2}{3}\eps\right\}\right|
> \frac{2}{3} (2^{\Delta h}),
\end{equation*}
and for all $j$ such that
\begin{equation*}
|\bar\weight(a_{j},b_{j}) 
- \weight(a_{0},b_{0})| > \eps,
\end{equation*}
we have that
\begin{equation*}
\left|\left\{j'\,:\,|\bar\weight(a_{j'},b_{j'}) 
- \bar\weight(a_{j},b_{j})| < \frac{2}{3}\eps\right\}\right|
\leq \frac{1}{3} (2^{\Delta h}).
\end{equation*}
Finally, we get that
\begin{equation*}
|\bar\weight'(a_0,b_0)
- \weight(a_{0},b_{0})| \leq \eps.
\end{equation*}
\end{app-proof}

\begin{app-proof}{Proposition~\ref{lem:deep2}}
The proof is similar to the proof
of Lemma~\ref{lem:distmet2gtr} and Proposition~\ref{lem:deep1}.
\end{app-proof}

\begin{app-proof}{Proposition~\ref{prop:deep3}}
We let $\eps < \quantum /2$.

The first part of the proposition
follows immediately from Proposition~\ref{lem:deep1} and the second part 
of Proposition~\ref{lem:deep2}.
Note that Propositions~\ref{lem:deep1} and~\ref{lem:deep2}
are still valid when $h' < \loglog{\alpha}$ if we 
take instead $\Delta h = \min\{h', \loglog{\alpha}\}$.
Indeed, in that case there is no need to perform an implicit
ancestral sequence reconstruction and the proofs of the lemmas
follow immediately from Lemmas~\ref{lem:distmet1gtr} and~\ref{lem:distmet2gtr}.

For the second part, choose $\kappa$ so as to satisfy the conditions of 
Proposition~\ref{lem:deep1}
{\em with diameter $D + \ln W$} and apply the first part of Proposition~\ref{lem:deep2}.
\end{app-proof}

\begin{app-proof}{Proposition~\ref{prop:weights}}
The proof follows immediately from Proposition~\ref{prop:deep3} and the remark above the statement of
Proposition~\ref{prop:weights}.
\end{app-proof}

\begin{app-proof}{Theorem~\ref{thm:mainhmg}} The proof of Theorem~\ref{thm:mainhmg} follows
from Propositions~\ref{prop:deep3} and~\ref{prop:weights}.
Indeed, at each level $h'$, 
we are guaranteed by the above to compute a distorted metric with a radius
large enough to detect all cherries on the next level using four-point tests.
The proof follows by induction. 
\end{app-proof}

\end{document}